\documentclass[12pt, a4paper, reqno, oneside]{amsart}      
\usepackage[left=2.20cm, right=2.2cm, top=2.0cm, bottom=1.8cm]{geometry}
\usepackage{amsthm}
\usepackage{bbm}
\usepackage{dsfont}

\usepackage{longtable}
\usepackage{mathrsfs}
\usepackage{mathtools}
\usepackage{amsmath,amssymb}
\newcommand\bigzero{\makebox(0,0){\text{\huge0}}}
\usepackage{amsthm}
\usepackage{bbold}
\usepackage{dsfont}
\usepackage{float}

\usepackage{graphicx}               
\usepackage{amssymb}

\newcommand{\R}{\mathbb{R}}

\newcommand{\eps}{\varepsilon}

\newcommand{\dd}{{\rm d}}

\usepackage{color}
\usepackage[unicode]{hyperref}
\hypersetup{
	colorlinks = true,%
	citecolor = [rgb]{0.0,0.0,0.9},
	filecolor=black,%
	linkcolor = [rgb]{0.65,0.0,0.0},%
	anchorcolor = red,
	pagecolor = red,
	urlcolor= [rgb]{0.65,0.0,0.0}
}

\usepackage{amssymb,amsmath}
\usepackage{graphicx}
\usepackage{mathrsfs}
\usepackage{float}
\usepackage[active]{srcltx}
\usepackage{caption}

\numberwithin{equation}{section}

\newtheorem{example}{ Example}[section]
\newtheorem{proposition}{Proposition}[section]
\newtheorem{theorem}{Theorem}[section]
\newtheorem{lemma}{Lemma}[section]
\newtheorem{corollary}{Corollary}[section]

\newtheorem{remark}{Remark}[section]
\newtheorem{claim}{Claim}[section]

\makeatletter
\def\l@section{\@tocline{1}{0pt}{1pc}{}{}}
\def\l@subsection{\@tocline{2}{0pt}{1pc}{4.6em}{}}
\def\l@subsubsection{\@tocline{3}{0pt}{1pc}{7.6em}{}}
\renewcommand{\tocsection}[3]{%
	\indentlabel{\@ifnotempty{#2}{\makebox[2.3em][l]{%
				\ignorespaces#1 #2.\hfill}}}#3}
\renewcommand{\tocsubsection}[3]{%
	\indentlabel{\@ifnotempty{#2}{\hspace*{2.3em}\makebox[2.3em][l]{%
				\ignorespaces#1 #2.\hfill}}}#3}
\renewcommand{\tocsubsubsection}[3]{%
	\indentlabel{\@ifnotempty{#2}{\hspace*{4.6em}\makebox[3em][l]{%
				\ignorespaces#1 #2.\hfill}}}#3}
\makeatother

\title[Jacobian determinant inequality]
{Jacobian determinant inequality on corank 1 Carnot groups with 
	applications}

\makeatletter
\renewcommand{\@biblabel}[1]{[#1]\hfill}
\makeatother

\author{Zolt\'an M. Balogh, Alexandru Krist\'aly, and Kinga Sipos}	

\thanks{Z. M. Balogh was
	supported by the Swiss National Science Foundation, Grant Nr. {200020\_165507}.  A. Krist\'aly  was supported by the STAR-UBB Advanced Fellowship-Intern, (Project PN II1.2PDI-PFC-C1-PFE-404).
	K. Sipos was supported by ERC Marie-Curie Research and Training Network MANET}

\begin{document}
\maketitle
\vspace*{-0.5cm}
\begin{abstract}
	We establish a weighted pointwise Jacobian determinant inequality on corank 1 Carnot groups related to optimal mass transportation akin to the work of Cordero-Erausquin, McCann and Schmuckenschl\"ager. In this setting, the presence of  abnormal geodesics does not allow the application of the general sub-Riemannian optimal mass transportation theory developed by  Figalli and Rifford and we need to work with a weaker notion  of Jacobian determinant. Nevertheless, our result achieves a transition between  Euclidean and sub-Riemannian structures, corresponding to the mass transportation along abnormal and strictly normal geodesics, respectively. The weights appearing in our expression are distortion coefficients that  reflect the delicate sub-Riemannian structure of our space.  As applications, entropy, Brunn-Minkowski and Borell-Brascamp-Lieb inequalities are established on Carnot groups. 
\end{abstract}

	\vspace*{0.5cm}




\section{Introduction}
As a general framework of our results, let  $(X,d,\textsf{m})$ be a suitably regular geodesic metric measure space with topological dimension $N\in \mathbb N$ where the theory of optimal mass transportation can be successfully developed.  Examples for such spaces include Riemannian and Finsler manifolds, see McCann \cite{McCann-GAFA} and  Ohta \cite{Ohta}, the Heisenberg group $\mathbb H^n$, see Ambrosio and Rigot \cite{AR}, or even more general sub-Riemannian structures with 'well-behaved' cut locus, see Figalli and Rifford \cite{FR}. 
Let $\mu_0$ and $\mu_1$ be two probability measures on $X$ which are absolutely continuous w.r.t. the reference measure $\textsf{m}$,  and let $\mu_s=(\psi_s)_\#\mu_0,$ $s\in [0,1],$ be the unique displacement interpolation measure joining  $\mu_0$ and $\mu_1$ throughout the so-called $s$-intermediate optimal transport map $\psi_s:X\to X$. Roughly speaking, for $s\in (0,1)$ fixed, the Jacobian determinant inequality reads as
\begin{equation}\label{Jacobian-000}
\left({\rm Jac}(\psi_s)(x)\right)^\frac{1}{N}\geq \tau_{1-s}^{N}(\theta_x)+\tau_{s}^{N}(\theta_x)\left({\rm Jac}(\psi)(x)\right)^\frac{1}{N}\mbox{ for } \mu_0 \mbox{-a.e. } x \in X.
\end{equation}
Here, {and in the sequel} ${\rm Jac}(\psi_s)(x)$ and ${\rm Jac}(\psi)(x)$ are interpreted as {densities, or} the Radon-Nikodym derivatives of $\mu_s$ and of $\mu_1$ w.r.t. the reference measure $\textsf{m}.$  {Note that in case when $X= \R^n$ and $\psi_s$ is differentiable at $x$ the term  ${\rm Jac}(\psi_s)(x)$ can be computed as ${\rm Jac}(\psi_s)(x)= |\det D \psi_s(x) |$. On the other hand, the Jacobian determinant in the above sense might exist as density even in the case when $\psi_s$ is not differentiable.}  The expression  $\tau_{s}^{N}$ is the distortion coefficient which encodes information on the geometric structure of the space $X$.  Expressions of $\tau_{s}^{N}$ can be calculated in terms of the Jacobian of the exponential map or estimated in terms of a curvature condition. The expression $\theta_x$ can be given as a function of $d(x,\psi(x))$ or its derivatives. 
 
The Jacobian determinant inequality (\ref{Jacobian-000}) in the above general form has been considered first in the setting of complete Riemannian manifolds (endowed with the natural Riemaniann distance and volume form) in the pioneering work of Cordero-Erausquin, McCann and Schmuckenschl{\"a}ger \cite{CMS}.  This result constituted the starting point of an extensive  study of the geometry of metric measure spaces, while relation (\ref{Jacobian-000}) became an equivalent formulation of the famous curvature-dimension condition $CD(K,N)$, due to  Lott and Villani \cite{LV}, and Sturm \cite{Sturm1, Sturm2}, where  $\tau_{s}^{N}$ is replaced by explicit expressions $\tau_{s}^{N, K}$, $K$ being the lower bound of the Ricci curvatures in the Riemannian setting. Namely, $\tau_{s}^{N, K}$ is given by
{\small
		$$\tau_s^{K,N}(\theta)=\left\{
		\begin{array}{lll}
	 s^\frac{1}{N}\left(\sinh\left(\sqrt{-\frac{K}{N-1}}s\theta\right)\big/\sinh\left(\sqrt{-\frac{K}{N-1}}\theta\right)\right)^{1-\frac{1}{N}}&
	 {\rm if} & K\theta^2<0;\\
	 	s & {\rm if} & K\theta^2=0;\\
	 s^\frac{1}{N}\left(\sin\left(\sqrt{\frac{K}{N-1}}s\theta\right)\big/\sin\left(\sqrt{\frac{K}{N-1}}\theta\right)\right)^{1-\frac{1}{N}}& {\rm if} &
		0<K\theta^2<(N-1)\pi^2;\\
			+\infty & {\rm if} & K\theta^2\geq (N-1)\pi^2,
		\end{array}
		\right.$$}
and $\theta=\theta_x$ is precisely the Riemannian distance $ d(x,\psi(x))$.

Juillet \cite{Jui} proved that the Lott-Sturm-Villani curvature-dimension condition does not hold for any  pair of parameters $(N,K)$ on the  Heisenberg group $\mathbb H^n$ (endowed  with its usual Carnot-Carath\'eo\-do\-ry metric $d_{CC}$ and $\mathcal L^{2n+1}$-measure), which is the simplest sub-Riemannian structure. 
Accordingly, there were strong doubts on the validity of a sub-Riemannian version of the Jacobian determinant inequality in the sub-Riemannian context.  However, by using a natural Riemannian approximation of the Heisenberg group as in Ambrosio and Rigot \cite{AR}, the authors of the present paper proved (\ref{Jacobian-000}) on $\mathbb H^n$, see \cite{BKS1,BKS}, where the Heisenberg distortion coefficient 
  $\tau_s^{2n+1}:[0,2\pi]\to [0,\infty]$ is defined by
  \begin{eqnarray}\label{concentration}
  \tau_s^{2n+1}(\theta) = \left\{
  \begin{array}{lll}
  {s^\frac{1}{2n+1}}
  \left(\frac{\sin\frac{\theta s}{2}}{\sin\frac{\theta }{2}}\right)^\frac{2n-1}{2n+1}\left(\frac{\sin\frac{\theta s}{2}-\frac{\theta s}{2}\cos\frac{\theta s}{2}}{\sin\frac{\theta }{2}-\frac{\theta }{2}\cos\frac{\theta
  	}{2}}\right)^\frac{1}{2n+1}
  \  &\mbox{if} &  \theta\in( 0,2\pi); \\
  s^\frac{2n+3}{2n+1} &\mbox{if} &  \theta=0;\\
   +\infty &\mbox{if} &  \theta=2\pi,\\
  \end{array}\right.
  \end{eqnarray}
    and $\theta = \theta_x$ is the 'vertical' derivative of $\frac{d_{CC}^2(\psi(x),\cdot)}{2}$ at the point $x$.

In the present paper we prove a Jacobian determinant inequality on corank 1 Carnot groups where the  sub-Riemannian geometry is more complicated  than the one of the model Heisenberg group $\mathbb H^n$  due to the presence of {abnormal} geodesics  and the 'anisotropic' structure of the cut locus.  
Our method is different from the one in \cite{BKS1, BKS} as we obtain the  Jacobian determinant inequality by an intrinsic approach, without using a Riemannian approximation.   As in \cite{BKS1, BKS}, we  apply our Jacobian determinant inequality to establish various functional and geometric inequalities in the present setting including entropy, Brunn-Minkowski and Borell-Brascamp-Lieb inequalities. These results should open up the way to considering the above inequalities in a broader context outside the realm of $CD(K,N)$-type conditions by replacing the coefficients $\tau_{s}^{N, K}$ by expressions that are suitable for sub-Riemannian geometries. {In this way, our results motivate the so-called "grande unification" of the three main geometries (Riemannian, Finslerian and sub-Riemannian), suggested by C. Villani in \cite[p. 43]{Villani-2}.}

In order to present our main result, let us fix some notation. We denote by $G$ a $k+1$ dimensional corank 1 Carnot group with its Lie algebra $\mathfrak g=\mathfrak g_1\oplus \mathfrak g_2$, where  dim$\mathfrak g_1=k\geq 2$ and dim$\mathfrak g_2=1$. The operation on $\mathfrak g$ (in exponential coordinates on $\mathbb R^k\times \mathbb R$) can be given by 
$$x\circ y=\left(x_1+y_1,...,x_k+y_k,x_z+y_z-\frac{1}{2}\sum_{i,j=1}^k \mathcal A_{ij}x_jy_i\right ),$$
where $x=(x_1,...,x_k,x_z)$, $y=(y_1,...,y_k,y_z)$, and  $\mathcal A=[\mathcal A_{ij}]$ is a $k\times k$ real skew-symmetric matrix. Let $e=(0_{\mathbb R^k},0)\in \mathbb R^k\times \mathbb R$ be the neutral element in $(G,\circ).$
The layers $\mathfrak g_1$  and $\mathfrak g_2$ are generated by the left-invariant vector fields 
\begin{equation}\label{vector-field}
X_i=\partial_{x_i}-\frac{1}{2}\sum_{j=1}^k \mathcal A_{ij}x_j\partial_{z},\ \ i=1,...,k.
\end{equation}
Moreover,  $[X_i,X_j]=\mathcal A_{ij}\partial_{ z}.$ 
By the spectral theorem for skew-symmetric matrices one can consider the diagonalized representation of $\mathcal A$  given by
\begin{equation}\label{matrix-representation}
\mathcal A=\left[
\begin{array}{cccc}
\mathbb{0}_{k-2d}& & &  \bigzero\\
& \alpha_1 J & & \\
\bigzero &   & \ddots \\
& & &  \alpha_d J
\end{array}
\right],\ \ \ J=\left[\begin{matrix}
0 & 1\\
-1& 0
\end{matrix}\right],
\end{equation}
where $0<\alpha_1\leq ...\leq \alpha_d,$ and $\mathbb{0}_{k-2d}$ is the $(k-2d)\times (k-2d)$ square null-matrix; {from now on, we assume the matrix $\mathcal A$ has this  representation. }

For further use, let us introduce the functions $\mathbb d_1, \mathbb d_2:[0,2\pi]\times(0,1)\to \mathbb R$ given by 
$$\mathbb d_1(t,s)=\frac{\sin(ts/2)}{s}\ \ {\rm and}\ \ \mathbb d_2(t,s)=\frac{\sin(ts/2)-ts/2\cos(ts/2)}{s}.$$ To define the distortion coefficient, we introduce the set
 $$D = \left\{p = (p_x,p_z) \in \mathbb{R}^{k+1} : |p_z| < \frac{2\pi}{\alpha_d} \mbox{ and } \mathcal Ap_x \neq 0_{\mathbb R^k}\right\} \subset T_e^* G,$$ 
 where $p_x=(p_x^0,p_x^1,...,p_x^d)\in \mathbb R^{k-2d}\times \mathbb R^2\times...\times \mathbb R^2$, and let 
 $\overline D$ be the closure of $D$.  The {\it distortion coefficient $\tau_s^{k,\alpha}:\overline D\to \mathbb R$ on the Carnot group} $(G,\circ)$ is defined by
{\small \begin{eqnarray*}
		\tau_s^{k,\alpha}(p)=  \left\{
		\begin{array}{lll}
			s\left(\frac{\displaystyle\sum_{i=1}^d\|p_x^i\|^2
				\displaystyle\prod_{j\neq i}\mathbb d_1^2(\alpha_j p_z,s)\mathbb d_1(\alpha_i p_z,s)\mathbb d_2(\alpha_i p_z,s)}{\displaystyle\sum_{i=1}^d\|p_x^i\|^2
				\displaystyle\prod_{j\neq i}\mathbb d_1^2(\alpha_j p_z,1)\mathbb d_1(\alpha_i p_z,1)\mathbb d_2(\alpha_i p_z,1)}\right)^\frac{1}{k+1}
			
			\  &\mbox{if} &  p\in D \ \& \ p_z\neq 0; \\
			s^\frac{k+3}{k+1} &\mbox{if} &  p\in D \ {\rm \&}\  p_z= 0;\\
			+\infty&\mbox{if} &  \mathcal  Ap_x\neq 0_{\mathbb R^k} \ {\rm \&}\  |p_z|=\frac{2\pi}{\alpha_d};\\
			s&\mbox{if} & \mathcal  Ap_x=0_{\mathbb R^k},
		\end{array}\right.
	\end{eqnarray*}
}
%
where $p=(p_x,p_z)$ and $\alpha=(\alpha_1,...,\alpha_d).$ {The functions  $\mathbb d_1$ and $\mathbb d_2$ appear explicitly in the Jacobian of the exponential map, see (\ref{Jacobian-Juillet}) below. In fact, $\mathbb d_2$ is a typical sub-Riemannian function appearing once after differentiating the exponential map along the 'vertical' direction, while $\mathbb d_1$ appears on the diagonal of the Jacobian matrix with multiplicity $2d-1$, see also Rizzi \cite{Rizzi}.
}

Let us consider two compactly supported probability measures $\mu_0$ and $\mu_1$ on $G$ which are absolutely continuous w.r.t. $\mathcal{L}^{k+1}$. Since the distribution $\Delta=\{X_1,...,X_k\}$ on the corank 1 Carnot group $G$ is  two-generating,  
there exists a unique map realizing the optimal transportation between the measures $\mu_0$ and $\mu_1$ w.r.t. the cost function ${d^2_{CC}}/{2}$, see Figalli and Rifford 
\cite[Proposition 4.2 and Theorem 3.2]{FR}; this map can be defined $\mu_0$-a.e. through 
a $d_{CC}^2/2$-concave function  $\varphi:G\to \mathbb R$ as 
{\begin{eqnarray} \label{DefIntMapH}
	\psi(x):=
	\left\{
	\begin{array}{lll}
	 \exp_x(-\nabla\varphi(x)) 
	\  &\mbox{if} &  x\in \mathcal M_\varphi\cap {\rm supp}(\mu_0); \\
	x &\mbox{if} &  x\in \mathcal S_\varphi\cap {\rm supp}(\mu_0).
	\end{array}\right.
	\end{eqnarray}
Hereafter, $d_{CC}$ is the Carnot-Carath\'eo\-do\-ry metric on $G$ and the sets $\mathcal M_\varphi$ and $\mathcal S_\varphi$ denote the moving and static sets of the transportation, respectively; see Section \ref{SecPrelim} for details. 
For $s\in (0,1)$ fixed, we also introduce the $s$-interpolant optimal transport map as
{\begin{eqnarray} \label{DefIntMapH-2}
	 \psi_s(x):=
	\left\{
	\begin{array}{lll}
	 \exp_x(-s\nabla\varphi(x)) 
	\  &\mbox{if} &  x\in \mathcal M_\varphi\cap {\rm supp}(\mu_0); \\
	x &\mbox{if} &  x\in \mathcal S_\varphi\cap {\rm supp}(\mu_0).
	\end{array}\right.
	\end{eqnarray}

Our main result  reads as follows.

\begin{theorem}\label{TJacobianDetIneq}{\bf (Jacobian determinant inequality on  Carnot groups)} Let $(G,\circ)$ be a $k+1$ dimensional corank 1 Carnot group,   and assume that $\mu_0$ and $\mu_1$ are two compactly supported Borel probability measures on  $G$, both absolutely continuous w.r.t. $\mathcal L^{k+1}$.  Let $s \in (0,1)$ be fixed, $\psi:G\to G$ be the unique optimal transport map transporting $\mu_0$ to $\mu_1$ associated to the cost function $\frac{d^2_{CC}}{2}$ and $\psi_s$ its $s$-interpolant map. Then the following Jacobian determinant inequality holds
	\begin{equation}\label{Jacobi-inequality-elso}
	\left({\rm Jac}(\psi_s)(x)\right)^\frac{1}{k+1}\geq \tau_{1-s}^{k,\alpha}(\theta_x)+\tau_{s}^{k,\alpha}(\theta_x)\left({\rm Jac}(\psi)(x)\right)^\frac{1}{k+1}\mbox{ for } \mu_0 \mbox{-a.e. } x \in G,
	\end{equation}
	where $\theta_x=(p_x,p_z)\in T_e^*G$ is given by $\exp_e(\theta_x)=x^{-1}\circ \psi(x).$
\end{theorem}
\medskip

Let us notice that if $p=(p_x,p_z)\in D$, we have that
 $$\lim_{p_z\to 0}\tau^{k,\alpha}_s(p)=s^{\frac{k+3}{k+1}}\ {\rm and}\ \lim_{p_z\to \pm2\pi/\alpha_d}\tau^{k,\alpha}_s(p)=+\infty.$$ Furthermore, monotonicity properties of the functions $\mathbb d_1$ and $\mathbb d_2$ {(cf. \cite[Lemma 2.1]{BKS})} show that 
 \begin{eqnarray}\label{tau-lower-bound}
 \tau^{k,\alpha}_s(p)\geq s^{\frac{k+3}{k+1}}\ {\rm for\ all}\ s\in (0,1),\ p\in \overline D.
 \end{eqnarray}
Therefore, the {\it measure contraction property} {\rm{\textsf{ MCP}}}$(0,k+3)$ proved by Rizzi \cite{Rizzi} is formally a consequence of (\ref{Jacobi-inequality-elso}). 
Notice, however that we use Rizzi's result to prove the absolute continuity of the interpolant measure $\mu_s=(\psi_s)_\#\mu_0$ (see Proposition \ref{interpolant-absolute-cont}), needed in the proof of the Jacobian determinant inequality.

 

In our next remark we consider the situation when $G=\mathbb H^n$ is the $n$-dimensional Heisenberg group. In this case we have $k=2n=2d$ and $\alpha_i=4$ for every $i\in \{1,...,d\}$. Moreover, no abnormal geodesics appear in $\mathbb H^n$  and the Carnot distortion coefficient $\tau_s^{2n,\alpha}(p_x,p_z)$ reduces to the Heisenberg distortion coefficient $\tau_s^{2n+1}(4p_z),$ which is nothing but relation \eqref{concentration} (introduced  in \cite{BKS}). Thus, most of the results of \cite{BKS} will be covered in the present work. 

Let us notice furthermore, that in general corank 1 Carnot groups, the coefficients $\tau^{k,\alpha}_s$ and $\tau^{k,\alpha}_{1-s}$  depend not only on the parameter $p_z$ (as in the Heisenberg group) but also on $\|p_x^i\|$, $i\in \{1,...,d\}$,  showing  a more anisotropic character of the present geometric setting as compared to the Heisenberg group. As we shall see later,  $\|p_x^i\|$ and $p_z$  can be obtained by differentiating $\frac{d_{CC}^2(\psi(x),\cdot)}{2}$ at the point $x$ w.r.t. the horizontal vector fields from the distribution $\Delta$ and the vertical vector field $\partial_z$, respectively (see Lemma \ref{LemmaFirstDeriv}  below).


Our final remark is of technical nature, but the details will be clear by reading the proof of Theorem \ref{TJacobianDetIneq}. In this proof,  we shall distinguish the cases when the mass is transported along \textit{{abnormal}} and \textit{{strictly normal}} geodesics, respectively. On  one hand, when the mass transport is realized along {\it abnormal} geodesics, it turns out that the Jacobian determinant inequality reduces to an {\it Euclidean-type} determinant inequality  thus the distortion coefficient can be  $\tau^{k,\alpha}_s=s$ as in the Euclidean framework. We notice that in this case the full Jacobian matrix of $\psi_s$ might not exist; however, since the matrix has a triangular structure, the Jacobian can be reduced to two parts of the diagonal which are well defined and inequality (\ref{Jacobi-inequality-elso}) makes sense. Furthermore, the triangular structure of the Jacobi matrix will allow us to perform the necessary changes of variable in order to provide important applications (see e.g. the entropy and Borell-Brascamp-Lieb inequalities via a suitable Monge-Amp\`ere equation). On the other hand,  once
the mass transport is along {\it strictly normal}
geodesics,  the distortion coefficient $\tau^{k,\alpha}_s$ encodes information on the genuine {\it sub-Riemannian} character of the Carnot group obtained by a careful analysis of the Jacobian for the exponential map. It could also happen that a positive part of the mass is transported along abnormal geodesics while the complementary mass is transported by strictly normal geodesics, so different formulas for  $\tau^{k,\alpha}_s$  will be used in the same instance of the mass transportation;  such a scenario will be presented in Example \ref{example} (see also Figure \ref{abra-elso}). {In conclusion, our results can be applied also in the presence of both abnormal and strictly normal geodesics in the so-called \textit{non-ideal} sub-Riemannian setting. Similar result in the case of general \textit{ideal} sub-Riemannian  geometries have been recently obtained by Barilari and Rizzi  \cite{Barilari-Rizzi}.}

The organization of the paper is as follows. The proof of Theorem \ref{TJacobianDetIneq} will be provided in Section \ref{SecProof} after a self-contained presentation of the needed technical details in Section \ref{SecPrelim}, i.e., properties of the Carnot-Carath\'eodory metric $d_{CC}$, exponential map and its Jacobian, the cut locus, and the optimal mass transportation on corank 1 Carnot groups. We emphasize that the optimal mass transportation developed by Figalli and Rifford \cite{FR} for large  classes of sub-Riemannian manifolds cannot be directly applied since the squared distance function $d_{CC}^2$ is not necessarily locally semiconcave outside of the diagonal of $G\times G$ which is crucial in \cite{FR} (e.g.  the regularity of  optimal mass transport maps $\psi$ and $\psi_s$, or the validity of the Monge-Amp\`ere equation). 
Section \ref{SecApps} is devoted to applications, i.e., by the Jacobian determinant inequality we shall derive entropy inequalities, the Brunn-Minkovski inequality and the Borell-Brascamp-Lieb inequality  on corank 1 Carnot groups. 

\medskip

\noindent {\bf Acknowledgements.} We express our gratitude to Luca Rizzi for motivating conversations about the subject of this paper. A. Krist\'aly is grateful to the Mathematisches Institute of Bern for the warm hospitality where this work has been developed. We also wish to thank the anonymous referees for their detailed reports and valuable
comments that greatly improved the presentation of the manuscript.

\section{Preliminaries} \label{SecPrelim}

\subsection{Carnot-Carath\'eodory metric and energy functional on {corank 1} Carnot groups} 


 {We shall consider a corank 1 Carnot group $(G,\circ)$, and make use of the notations already introduced in the previous section.} A horizontal curve on $(G,\circ)$ is an absolutely continuous curve $\gamma:
[0, r] \to G$ for which there exist {bounded} measurable functions ${u}_j:
[0,r] \to \R$ ($j = 1, ..., k$) such that
\begin{equation} \label{horiz}
\dot{\gamma}(s) = \sum\limits_{j=1}^k {u_j}(s) X_j(\gamma(s)) \quad \mbox{a.e. } s \in [0,r].
\end{equation}
In the sequel we denote by $\gamma_u$  such a horizontal curve. 
The length of this curve is given by 
$$l(u) = l(\gamma_u) = \int\limits_{0}^{r} \|\dot{\gamma_u}(s)\| \dd s= \int\limits_{0}^{r} \sqrt{\sum\limits_{j = 1}^{k} {u^2_j}(s) } \dd s.$$
The classical Chow-Rashewsky theorem assures that any two points
from the Carnot group can be joined by a horizontal curve. Thus
we can equip the Carnot group $G$ with its natural 
Carnot-Carath\'eodory metric by 
$$d_{CC}(x,y) = \inf \{l(\gamma): \gamma \mbox{ is a horizontal curve joining } x \mbox{ and } y\},$$
where  $x,y\in G$ are arbitrarily fixed.  

Let $e=(0_{\mathbb R^k},0) \in \mathbb{R}^k \times \mathbb{R}$ be the neutral element in $(G,\circ)$. 
The left invariance of the vector fields in the distribution $\Delta=\{X_1,
\dots, X_k\}$ is inherited by the distance $d_{CC}$,
thus
$$d_{CC}(x, y) = d_{CC}(e, x^{-1} \circ y) \quad {\rm for\ every}\ x,y \in G.$$
{Beside the length function $u\mapsto l(\gamma_u)$ we also consider the energy functional 
$$J(u) =\frac{1}{2} \int\limits_{0}^{r} \|\dot{\gamma_u}(s)\|^2 \dd s= \int\limits_{0}^{r} {\sum\limits_{j = 1}^{k} {u^2_j}(s) } \dd s.$$ It is well-known that the minimisers of $J$ induce {up to a reparametrisation} length minimising horizontal curves with constant speed between two  fixed endpoints.}

\subsection{Geodesics, exponential map and its Jacobian}

{Geodesics are horizontal curves that are locally energy minimizers between their endpoints. Let $\mathcal U\subset L^\infty([0,r],\mathbb R^k)$ be an open set and for a fixed $x\in G$, let  $E_x:\mathcal U\to G$ be the usual end-point map, $E_x(u)=\gamma_u(r)$, where $\gamma_u$ is the unique curve with the property that $\gamma_u(0)= x$ and satisfying \eqref{horiz} see e.g. Figalli and Rifford \cite[\S 2.1]{FR}. A minimizing geodesic $\gamma_u $ for $ u\in \mathcal U$ is a solution of the problem
$$J(v)\to\min,\ \ E_x(v)=y,\ \ v\in \mathcal U.$$
According to the Lagrange multipliers rule, there is $(\lambda,\mu)\in T^*_yG\times \{0,1\}\setminus \{ (0,0)\}$ such that 
$$\lambda (D_uE_x)=\mu D_uJ.$$ The associated curve $\gamma_u$ is normal if $\mu=1$ and abnormal if $\mu=0$ (the latter being equivalent to the fact that $u$ is a critical point of $E_x$). We notice that on any corank 1 Carnot group all minimizing geodesics are normal.  
Following Rizzi \cite{Rizzi}, the explicit form of such {normal} minimal geodesics  can be described as follows.}

 
\begin{proposition} {\rm (Rizzi \cite{Rizzi})}\label{proposition-geodetikus}
On a corank $1$ Carnot group $(G,\circ)$  the geodesic $s\mapsto \exp_e(sp)\in G$  starting from $e=(0_{\mathbb R^k},0)$, with {initial covector}
$$p=(\underbrace{p_x^0, p_x^1, ..., p_x^d}_{p_x},p_z) \in \left(\mathbb{R}^{k-2d} \times \mathbb{R}^2 \times ... \times \mathbb{R}^2\right) \times \mathbb{R}= T_e^*G $$ 
has the following equation
\begin{eqnarray}\label{GeodesicEqOnCarnotGroup}
\exp_e(sp):
\left\{
\begin{array}{l}
\gamma^0(s) = p_x^0 s,\\
\gamma^i(s) = \left(\frac{\sin(\alpha_ip_z s)}{\alpha_i p_z}I + \frac{\cos(\alpha_i p_z s) -1}{\alpha_i p_z}J\right) p_x^i,\\
\gamma_z(s) = \sum_{i=1}^d \|p_x^i\|^2 \frac{\alpha_i p_z s - \sin(\alpha_i p_z s)}{2 \alpha_i p_z^2},
\end{array}
\right. \quad  s \in [0,1],
\end{eqnarray}
when $p_z\neq 0$. When $p_z = 0$, the geodesic is 
\begin{eqnarray}\label{GeodesicEqNullpz}
s\mapsto \exp_e(sp)=(p_x^0s, p_x^1s, ..., p_x^d s,0), \quad s \in [0,1].
\end{eqnarray}
Hereafter, $I$ denotes the $2\times 2$ unit matrix { and $J=\left[\begin{matrix}
	0 & 1\\
	-1& 0
\end{matrix}\right]. $}
\end{proposition}

{Once $\mathcal A$ has a non-trivial kernel, every nonzero covector $(p_x,p_z)$ with  $\mathcal Ap_x=0$, corresponds to an abnormal geodesic; more precisely, for every choice of $p_z\in \mathbb R$ one has 
\begin{eqnarray}\label{AbnGeodesicEq}
s\mapsto \exp_e(p_x^0s,0,...,0,p_zs)=(p_x^0s,0_{\mathbb R^{2d+1}}),\ s\in [0,1].
\end{eqnarray}
Note that the image of such a geodesic can be also obtained by (\ref{GeodesicEqNullpz}), letting   $p_z = 0$ and $p_x^1 = ... = p_x^d = 0_{\mathbb R^2}$. {These type of geodesics are normal and also abnormal at the same time}. It turns out that all abnormal geodesics have this representation.}

We recall from Rizzi \cite{Rizzi} that the Jacobian determinant of the exponential map is
\begin{eqnarray}\label{Jacobian-Juillet}
{\rm Jac}(\exp_e)(p) = \left\{
\begin{array}{lll}
\displaystyle\frac{2^{2d}}{\prod_{i=1}^d \alpha_i^2 p_z^{2d+2}} \displaystyle\sum_{i=1}^d\|p_x^i\|^2
\displaystyle\prod_{j\neq i}\left({\sin\frac{\alpha_j p_z }{2}}\right)^2\sin\frac{\alpha_i p_z }{2}\times\\
\ \ \ \ \ \ \ \ \ \ \  \ \ \ \ \ \ \ \ \ \ \ \ \ \ \ \ \ \ \ \ \ \  \times\left({\sin\frac{\alpha_i p_z }{2}-\frac{\alpha_i p_z }{2}\cos\frac{\alpha_i p_z }{2}}\right)
\  &\mbox{if} &  p_z\neq 0; \\
& & \\
\frac{1}{12}{\displaystyle\sum_{i=1}^d\|p_x^i\|^2}\alpha_i^2 &\mbox{if} &  p_z=0.
\end{array}\right.
\end{eqnarray}

%

\noindent 
By left-invariance, the minimal geodesics on $G$ starting from an arbitrary point $x\in G$ are represented by $s\mapsto \exp_x(sp)=x\circ \exp_e(s\tilde p)$, $s\in [0,1],$ where the two covectors $p\in T_x^*G$ and  $\tilde p\in T_e^*G$ can be identified. 
Moreover, since for every $x\in G$ the left-translation $L_x(y)=x\circ y$, $y\in G,$ is a volume-preserving map, it follows that 
\begin{equation}\label{Jacobi-left}
{\rm Jac}(\exp_x)(p)={\rm Jac}(\exp_e)(p)\ {\rm for\ every}\ p\in T_x^*G.
\end{equation}

Given $x,y\in G$ and assume that $x=\exp_y(p)$ for some $p=(p_x,p_z)=({p_x^0, p_x^1, ..., p_x^d},p_z)\in T_y^*G$. Then $y=\exp_x(\overline p)$, where $\overline p=(\overline p^0_x,\overline p^1_x,...,\overline p^d_x,\overline p_z)$ is given by 
\begin{equation}\label{reverse-repres}
\left\{
\begin{array}{lll}
\overline p^0_x=- p_x^0;\\
\overline p^i_x=\left(-\cos(\alpha_i p_z)I+\sin(\alpha_ip_z)J\right) p_x^i,\ \ i\in \{1,...,d\};\\
\overline p_z=- p_z.
\end{array}
\right.
\end{equation}

 We notice that  $\Delta=\{X_1,...,X_{k}\}$ is {\it not} a fat distribution whenever the kernel of $\mathcal A$ is non-trivial. Indeed, in this case we have $T_xG\neq \Delta(x)+[X_j,\Delta](x)$ for every $x\in G$ and $j\in \{1,...,k-2d\}$. However, $\Delta$ is two-generating, i.e.,  $$T_xG= \Delta(x)+[\Delta,\Delta](x)\ {\rm for\ every}\ x\in G.$$  
 For simplicity of notation, we reorganize the vector fields in $T_xG$ as
 \begin{equation}\label{vector-fields}
 \left\{
 \begin{array}{lll}
 X^0=(X_1,..., X_{k-2d});\\
 X^i=(X_{k-2d+2i-1}, X_{k-2d+2i}),\ \ i\in \{1,...,d\};\\
 Z=\partial_z.
 \end{array}
 \right.
 \end{equation}
 We  split the distribution $\Delta$ on $G$ into two types of vector fields; namely, $\Delta_0=\{X^0\}$ and $\tilde \Delta=\{X^{1},...,X^{d}\}$. This splitting gives the following trivial representation of the distance function $d_{CC}$:
 
 \begin{lemma}{\rm (Pythagorean rule)} \label{LemmaPythagorean}
  For every $(\xi,\eta,z), (\overline \xi,\overline \eta,\overline z)\in \mathbb R^{k-2d}\times \mathbb R^{2d}\times \mathbb R$ , we have 
 $$d_{CC}^2((\xi,\eta,z), (\overline \xi,\overline \eta,\overline z))=d^2_{\mathbb R^{k-2d}}(\xi,\overline \xi)+\tilde d_{CC}^2((\eta,z), (\overline \eta,\overline z)),$$ where $d_{\mathbb R^{k-2d}}$ is the Euclidean metric in $\mathbb R^{k-2d}$ while $\tilde d_{CC} $ is the Carnot-Carath\'eodory distance on $\mathbb R^{2d}\times \mathbb R$ w.r.t. to the distribution $\tilde \Delta$ inherited from the original sub-Riemannian structure.
 \end{lemma}
{\it Proof.} 
  By the left-invariance of the metric $d_{CC}$, we have 
  
  $$d_{CC}^2((\xi,\eta,z), (\overline \xi,\overline \eta,\overline z))=d_{CC}^2(e, (-\xi, -\eta, -z) \circ (\overline \xi, \overline \eta,\overline z)).$$
 Let $\gamma=(\gamma^0,\gamma^1,...,\gamma^d,\gamma_z):[0,1]\to G$ be the geodesic given by (\ref{GeodesicEqOnCarnotGroup}) or (\ref{GeodesicEqNullpz}) joining $e$ and the element $(-\xi, -\eta, -z) \circ (\overline \xi, \overline \eta,\overline z)$,  having its initial vector $p=(p_x^0,p_x^1,...,p_x^d,p_z)\in \mathbb R^{k-2d}\times \mathbb R^{2}\times \cdots \times \mathbb{R}^2 \mathbb \times \mathbb R$.  We have that $d_{CC}^2((\xi,\eta,z), (\overline \xi,\overline \eta,\overline z))=\sum_{i=0}^d \|p_x^i\|^2.$
 Note that $\|p_x^0\|_{\mathbb R^{k-2d}}=d_{\mathbb R^{k-2d}}(\xi,\overline \xi)$ and $$\sum_{i=1}^d \|p_x^i\|^2= d_{CC}^2(e, (0_{\mathbb R^{k-2d}}, -\eta, -z) \circ (0_{\mathbb R^{k-2d}}, \overline \eta,\overline z))=\tilde d_{CC}^2( (\eta, z),  (\overline \eta,\overline z))$$ which is realized precisely by the geodesic $\tilde \gamma=(\gamma^1,...,\gamma^d,\gamma_z)$, concluding the proof.
 \hfill $\square$

 \subsection{Cut locus} \label{SubseqCut}
 Let us consider the set
 $$D = \left\{p = (p_x,p_z) \in \mathbb{R}^{k+1} : |p_z| < \frac{2\pi}{\alpha_d} \mbox{ and } \mathcal  Ap_x \neq 0_{\mathbb R^k}\right\} \subset T_e^* G.$$ 
 Rizzi \cite[Lemma 16]{Rizzi} proved that $D$ is precisely the injectivity domain of parameters associated to  geodesics joining the origin $e$ {to almost all points of} $G$. We know that all points in the corank 1 Carnot group $G$ can be reached by a minimal normal geodesic; namely, for every $x \in G$ there exists a parametrization $p$ in the closure of $D$, i.e., 
 $$\overline{D} = \left\{p = (p_x,p_z) \in \mathbb{R}^{k+1} : |p_z| \leq \frac{2\pi}{\alpha_d}\right\},$$ 
 which defines a minimal normal geodesic joining $e$ and $x$. 
 
 The {\it cut locus} of the origin $e$ in $G$ is
 \begin{eqnarray*}
  {\rm cut}_G(e) &=& \exp_e(\overline{D} \setminus D)=G\setminus \exp_e(D)\\
 &=&\left(\mathbb R^{k-2d}\times \{0_{\mathbb R^{2d+1}}\}\right)\cup \left\{\exp_e\left(p_x,\pm\frac{2\pi}{\alpha_d}\right):\mathcal  Ap_x\neq 0_{\mathbb R^k}\right\}.
 \end{eqnarray*}
 The set $\mathbb R^{k-2d}\times \{0_{\mathbb R^{2d+1}}\}$ in the above representation corresponds to the image of abnormal geodesics while the latter set contains the conjugate points to $e$, see  (\ref{Jacobian-Juillet}).  
{Corank 1 Carnot groups have negligible cut loci, see Rizzi \cite[Section 1.4]{Rizzi}; alternatively, due to (\ref{GeodesicEqOnCarnotGroup}), one has that  ${\rm cut}_G(e)\subset \mathbb R^{k-2}\times 0_{\mathbb R^2}\times \mathbb R$ , thus $\mathcal L^{k+1}({\rm cut}_G(e))=0$.}
By left-invariance,  the cut locus of the point  $x\in G$ is $${\rm cut}_G(x) = L_x( {\rm cut}_G(e)),$$ 
 %
thus ${\rm cut}_G(x)$ is closed and $\mathcal L^{k+1}({\rm cut}_G(x))=0$ for every $x\in G;$ moreover, by (\ref{reverse-repres}) it follows that $y\in {\rm cut}_G(x)$ if and only if $x\in {\rm cut}_G(y)$.

 {The following two results are specifications to the case of corank 1 Carnot groups of the well-known fact $f:=\frac{1}{2}d_{CC}^2(y,\cdot)$ is smooth in a neighborhood
 of $x\in G$ whenever $x\notin {\rm cut}_G(y)$, and one can recover the initial covector $\lambda_x\in  T_xG$ of the unique
 geodesic joining $x$ with $y$ by  $\lambda_x=-\nabla f(x)$. }
 
\begin{lemma}\label{LemmaFirstDeriv} 
%
 
 Fix $y\in G$  and let $x=(x^0,x^1,...,x^d,z)\notin {\rm cut}_G(y)$. If $x=\exp_y(p_x^0,p_x^1,...,p_x^d,p_z)$ then  we have  
 \begin{itemize}
 	\item[(i)] $X^0\frac{d_{CC}^2(y,\cdot)}{2}\big|_{x}=p_x^0 $ and $Z \frac{d_{CC}^2(y,\cdot)}{2}\big|_{x}=p_z;$
 	\item[(ii)]  for every $i\in \{1,...,d\}$,  
\begin{equation}\label{derivalt-1}
 X^{i} \frac{d_{CC}^2(y,\cdot)}{2}\big|_x=[\cos(\alpha_i p_z)I-\sin(\alpha_i p_z)J]p_x^i.
 \end{equation}
\end{itemize}
 
 \end{lemma}
 
 {\it Proof.} By exploring the left-invariance, it is enough to consider the case when $y=e$.  Let us introduce the auxiliary functions $f,g:(-2\pi,2\pi)\setminus\{0\}\to \mathbb R$ defined by 
 \begin{equation}\label{fesg}
 f(t)=\frac{\sin^2\left(\frac{t}{2}\right)}{\left(\frac{t}{2}\right)^2}\ \ {\rm and}\ \ g(t)=\frac{t-\sin(t)}{\sin^2\left(\frac{t}{2}\right)},\ t\in (-2\pi,2\pi) \setminus\{0\}.
 \end{equation}
 
 We consider the case when $p_z\neq 0;$ the case $p_z=0$ can be obtained by a limiting procedure, {i.e., one must consider the limit $p_z\to 0$.} Since $x\notin {\rm cut}_G(e)$ and the cut locus is closed, there exists a small neighborhood $V_x$ of $x$ such that $V_x\cap {\rm cut}_G(e)=\emptyset$. Let $w=({x_w^0,x_w^1,...,x_w^d},z_w)={\exp_e\left((p_w)^0_x,(p_w)^1_x,...,(p_w)^d_x,(p_w)_z\right)}\in V_x$ be arbitrarily fixed.  By  (\ref{GeodesicEqOnCarnotGroup}) (for $s=1$) we have that 
 $$\|{x_w^i}\|^2={\|{(p_w)^i_x}\|^2}{f(\alpha_i {(p_w)_z})},\ i\in \{1,...,d\}.$$
 Thus, one has
 \begin{equation}\label{elso-dcc}
 d_{CC}^2(e,w)=\sum_{i=0}^d\|{(p_w)^i_x}\|^2=\|{x_w^0}\|^2+\sum_{i=1}^d\frac{\|{x_w^i}\|^2}{ f(\alpha_i {(p_w)_z})}.
 \end{equation}
 
(i)  By (\ref{elso-dcc}) we directly  have that $X^0(d_{CC}^2(e,\cdot))\big|_x=2x^0$. Furthermore, the last component in (\ref{GeodesicEqOnCarnotGroup}) can be written as
 \begin{equation}\label{masodik-dcc}
 z_w=\sum_{i=1}^d\|{(p_w)^i_x}\|^2\frac{\alpha_i {(p_w)_z}-\sin(\alpha_i {(p_w)_z})}{2\alpha_i {\left((p_w)_z\right)}^2}=\frac{1}{8}\sum_{i=1}^d\alpha_i\|{x_w^i}\|^2g(\alpha_i{(p_w)_z}).
 \end{equation}
 We may differentiate (\ref{elso-dcc}) and (\ref{masodik-dcc}) w.r.t. the variable  $z_w$ at  the point $x$, obtaining  $$Z(d_{CC}^2(e,\cdot))\big|_x=-\sum_{i=1}^d\alpha_i\|x^i\|^2 \frac{f'(\alpha_i p_z)}{f^2(\alpha_i p_z)}{\left(Z{(p_w)_z}\big|_x \right)} \ \ {\rm and} \ \ 1=\frac{1}{8}\sum_{i=1}^d\alpha_i^2\|x^i\|^2g'(\alpha_ip_z){\left(Z{(p_w)_z}\big|_x\right)}.$$
 Note that $-\frac{f'(t)}{f^2(t)}=\frac{t}{4}g'(t)$; thus, the latter relations give at once that
 $Z(d_{CC}^2(e,\cdot))\big|_x=2p_z.$

 (ii) In order to prove relation (\ref{derivalt-1})  we proceed in a similar way as in (i), by deriving (\ref{elso-dcc}) and (\ref{masodik-dcc}) w.r.t. the corresponding variables. \hfill$\square$\\

 A direct consequence of Lemma \ref{LemmaFirstDeriv} is: 
 
 \begin{proposition}\label{prop-carnot-exp}
 	Fix $x,y\in G$ such that $y \notin {\rm cut}_G(x)$. If $\nabla=({X^0,X^{1},...,X^{d}}, Z)$, then   
 	\begin{equation}\label{eqn-exp}
 	y=\exp_x\left(-\nabla\frac{d_{CC}^2(y,\cdot)}{2}\big|_x\right).
 	\end{equation}
 \end{proposition}
 
 {\it Proof.} Let  $x=\exp_y(p)$ for some $p=(p_x,p_z)=({p_x^0, p_x^1, ..., p_x^d},p_z)\in D$. According to Lemma \ref{LemmaFirstDeriv}, we have that 
 $$-\nabla\frac{d_{CC}^2(y,\cdot)}{2}\big|_x=(\overline p_x^0,\overline p_x^1,...,\overline p_x^d,\overline p_z),$$
 where 
$$ \left\{
 \begin{array}{lll}
 	\overline p_x^0=- p_x^0;\\
 	\overline p_x^i=-[\cos(\alpha_i p_z)I-\sin(\alpha_ip_z)J] p_x^i,\ \ i\in \{1,...,d\};\\
 	\overline p_z=- p_z.
 \end{array}
 \right.
 $$
 Thus, by relation (\ref{reverse-repres}) it follows that 
 $$\exp_x\left(-\nabla\frac{d_{CC}^2(y,\cdot)}{2}\big|_x\right)=\exp_x(\overline p_x^0,\overline p_x^1,...,\overline p_x^d,\overline p_z)=y,$$ 
 which concludes the proof.  \hfill$\square$\\

\subsection{The Jacobian of the exponential map along a reversed geodesic.}

Let $x,y\in G$ be such that $x\notin {\rm cut}_G(y)$ and $\gamma:[0,1]\to G$ be the unique geodesic $\gamma(s)=\exp_x(sp)$  joining $x$ and $y$ for some $p\in D$.   For every $s\in (0,1]$,  let us introduce the Jacobian matrix $$Y(s)=d(\exp_x)_{sp}.$$ According to (\ref{Jacobian-Juillet}), the matrix $Y(s)$ is invertible for  every $s\in (0,1]$. 
In the sequel, we are going to consider the reversed geodesic path $s\mapsto \exp_y((1-s) \overline p),$ $s\in [0,1],$ where $\exp_y \overline p = x$ and compute the 'reverse' of $Y$, i.e., \begin{equation}\label{tildeY}
\overline Y(1-s)=d(\exp_y)_{(1-s)\overline p}, \ s\in [0,1).
\end{equation}
Here, $\overline p\in T_y^*G$ is given by  $p\in T_x^*G$ similarly as in (\ref{reverse-repres}). With these  notations, we have 

\begin{proposition}\label{prop-hessian}
Let $x,y\in G$ be such that $x\notin {\rm cut}_G(y)$ and $\gamma:[0,1]\to G$ be the unique geodesic $\gamma(s)=\exp_x(sp)$  joining $x$ and $y$ for some $p\in T_x^*G$.	 For every $s\in (0,1),$ one has 	\begin{equation}\label{Y-tilde-bevezetese}
	 \overline Y(1-s)=\frac{1}{1-s}Y(s)H_{x,y}(s)\overline Y(1),
	 \end{equation} 
	 where
	 $$H_{x,y}(s)={\rm Hess}\frac{d_{CC}^2(\gamma(s),\cdot)}{2}\big|_{x}-s{\rm Hess}\frac{d_{CC}^2(y,\cdot)}{2}\big|_{x}.$$ In addition, $H_{x,y}(s)$ is a positive semidefinite, symmetric matrix.

\end{proposition}

  Let us note that in the above statement ${\rm Hess}=\nabla^2$ denotes the $($a priori not necessarily symmetric$)$ Carnot Hessian, i.e., 
	 $${\rm Hess}=\left[
	 \begin{array}{ccccc}
	 	X_1X_1& X_1X_2& ... &  X_1X_k& X_1Z\\
	 	X_2X_1& X_2X_2& ... &  X_2X_k& X_2Z\\
	 	\vdots& \vdots& ... &  \vdots& \vdots\\
	 	X_kX_1& X_kX_2& ... &  X_kX_k& X_kZ\\
	 	ZX_1& ZX_2& ... &  ZX_k& ZZ
	 \end{array}
	 \right].$$
	 This notation will be used also later on.

A similar result to Proposition \ref{prop-hessian} has been proved by Cordero-Erausquin,  McCann and Schmuckenschl{\"a}ger \cite{CMS} on Riemannian manifolds by exploring properties of Jacobi fields. Since the theory of Jacobi fields in our setting  is not (yet)  available, we give a direct proof of Proposition \ref{prop-hessian}. To do this, we need the following: 

\begin{claim}\label{claim-diff}
	Let $m\in \mathbb N$, $c, \eta_i:[0,1]\to \mathbb R^m$, $i\in \{1,2\},$ be some differentiable maps with $\eta_2(0)=0$ and a smooth function $F:\mathbb R^{2m}\to \mathbb R^m$ in a neighborhood of 
	$(c(0),\eta_1(0))$ such that $t\mapsto F(c(t),\eta_1(t))$ is constant near the origin. Then $$\frac{d}{dt}F(c(t),\eta_1(t)+\eta_2(t))|_{t=0}=D_2F(c(0),\eta_1(0)) \dot\eta_2(0).$$
\end{claim}

{\it Proof.} By assumption, we have near the origin that 
$$0=\frac{d}{dt}F(c(t),\eta_1(t))=D_1F(c(t),\eta_1(t)) \dot c(t)+D_2F(c(t),\eta_1(t))\dot \eta_1(t).$$  By using the latter relation at $t=0$ and $\eta_2(0)=0$, we obtain 
\begin{eqnarray*}
	\frac{d}{dt}F(c(t),\eta_1(t)+\eta_2(t))|_{t=0}&=&D_1F(c(0),\eta_1(0)) \dot c(0)+D_2F(c(0),\eta_1(0)) (\dot \eta_1(0)+\dot \eta_2(0))\\&=&D_2F(c(0),\eta_1(0)) \dot\eta_2(0),
\end{eqnarray*}
which completes the proof. 
\hfill$\square$\\


{\it Proof of Proposition \ref{prop-hessian}.}
We first deal with the properties of the matrix $H_{x,y}(s)$. By pure metric arguments, one can check that for every $z\in G$ and $s\in [0,1]$ we have the    inequality
\begin{eqnarray}\label{metric-MCS}
m_{x,y}^s(z):=d_{CC}^2(\gamma(s),z)/2- sd_{CC}^2(y,z)/2+s(1-s)d_{CC}^2(x,y)/2 \geq 0.
\end{eqnarray}
In the  Riemannian  setting this has been  established first
 by Cordero-Erausquin, McCann and Schmuc\-kenschl{\"a}ger \cite[Claim 2.4]{CMS}.  Moreover,  in (\ref{metric-MCS}) equality is realized precisely when $z = x$; the same proof works in our setting as well. 
 
  Since $\gamma((0,1])\cap {\rm cut}_G(x)=\emptyset,$ it follows that $z\mapsto m_{x,y}^s(z)$ is twice differentiable at $x$ (see Proposition \ref{prop-carnot-exp}) and its gradient is
\begin{equation}\label{derivative-null}
\nabla m_{x,y}^s(\cdot)|_x=\nabla\frac{d_{CC}^2(\gamma(s),\cdot)}{2}\big|_{x}-s\nabla\frac{d_{CC}^2(y,\cdot)}{2}\big|_{x}=0_{\mathbb R^{k+1}},
\end{equation}
 while its Carnot Hessian $\nabla^2 m_{x,y}^s(\cdot)|_x=H_{x,y}(s)$ is positive semidefinite.


In order to prove the symmetry of $H_{x,y}(s)$, we  verify that the Lie brackets $[W_1,W_2]m_{x,y}^s(\cdot)|_x$ vanish for every choice of $W_1,W_2\in \Delta\cup \{Z\}=\{X_1,...,X_k,Z\}.$ Indeed, the Lie bracket is either trivial by definition or it is  $Zm_{x,y}^s(\cdot)|_x$ up to a multiplicative constant (depending on the eigenvalues $\alpha_i$, $i\in \{1,...,d\}$); but $Zm_{x,y}^s(\cdot)|_x=0$ due to (\ref{derivative-null}). 

We now prove relation (\ref{Y-tilde-bevezetese}). Since $x\notin {\rm cut}_G(y)$ and  ${\rm cut}_G(y)$ is closed, one may fix a curve $c:[0,1]\to G$  with $c(0)=x$ and $\dot c(0)=w\in T_xG$   arbitrarily fixed such that $c([0,1])\cap {\rm cut}_G(y)=\emptyset$. 
We notice that $s\mapsto \exp_{c(t)}\left(-s\nabla\frac{d_{CC}^2(y,\cdot)}{2}\big|_{c(t)}\right){=:\gamma(s)}$ is the unique {minimal} geodesic joining $c(t)$ and $y$; indeed, for $s=0$ we have $c(t),$ while for $s=1$ one has precisely $y$ due to  Proposition \ref{prop-carnot-exp}, see Figure 
\ref{fig:YInvert}. {Moreover, by construction, it turns out that $\gamma(s)\notin {\rm cut}_G(c(t))$ for every $t,s\in[0,1]$.}

\begin{figure}[H]
 \centering
       \includegraphics[scale=0.4]{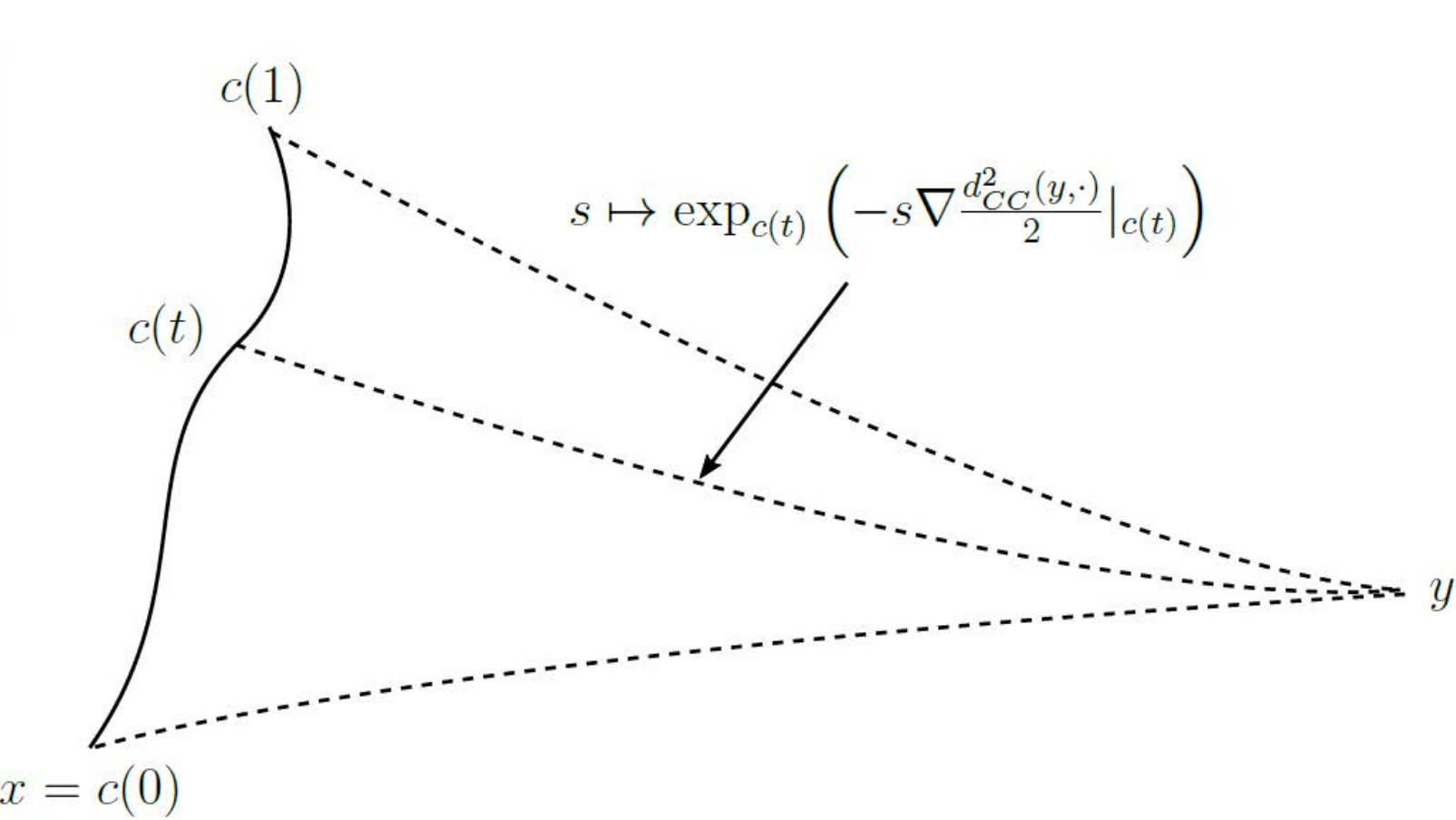}
        \caption{The curve $c$ (starting from $x$), connected by geodesics with the point $y$}
        \label{fig:YInvert}
\end{figure}

\noindent  Let 
$\overline p:[0,1]\to T_yG$ be a curve such that 
\begin{equation}\label{ketfele}
\exp_{c(t)}\left(-s\nabla\frac{d_{CC}^2(y,\cdot)}{2}\big|_{c(t)}\right)=\exp_y((1-s)\overline p(t)).
\end{equation}
Let us observe that by (\ref{ketfele}) for $s=0$ we have 
$$
c(t)=\exp_y(\overline p(t)) \ {\rm for\ all}\ t\in [0,1].
$$
In particular, for $t=0$  we have that  
$x=c(0)=\exp_y(\overline p(0))$, i.e., $\overline p(0)= \overline p$ and due to (\ref{tildeY}), 
\begin{equation}\label{wtildeY}
w=\dot c(0)=d(\exp_y)_{\overline p(0)}\dot{\overline p}(0)=d(\exp_y)_{\overline p}\dot{\overline p}(0)=\overline  Y(1)\dot{\overline p}(0).
\end{equation}
Fix $s \in (0,1)$. Now, we rewrite (\ref{ketfele}) into 
\begin{equation}\label{atirva}
\exp_{c(t)}\left(\eta_1(t)+\eta_2(t)\right)=\exp_y((1-s)\overline p(t)),
\end{equation}
where $$\eta_1^s(t)=-\nabla\frac{d_{CC}^2(\gamma(s),\cdot)}{2}\big|_{c(t)}\ {\rm and}\ \eta_2^s(t)=\nabla\frac{d_{CC}^2(\gamma(s),\cdot)}{2}\big|_{c(t)}-s\nabla\frac{d_{CC}^2(y,\cdot)}{2}\big|_{c(t)}.$$ We are going to verify the assumptions of Claim \ref{claim-diff} for the latter choices. 
First, due to Proposition \ref{prop-carnot-exp}, one has  $t\mapsto\exp_{c(t)}(\eta_1^s(t))=\gamma(s)=$constant, and due to (\ref{derivative-null}), we also have $\eta_2^s(0)=0.$ Since we have  $\eta_1^s(0)=-\nabla\frac{d_{CC}^2(\gamma(s),\cdot)}{2}\big|_{x}$,  by 
Proposition \ref{prop-carnot-exp} one has that $\exp_{x}(\eta_1^s(0))=\gamma(s)$ which is nothing but $\gamma(s)=\exp_x(sp)$; thus $\eta_1^s(0)=sp.$  Consequently, by  differentiating  relation (\ref{atirva}) at $t=0$ and using Claim \ref{claim-diff} with  $F(q_1,q_2)=\exp_{q_1}(q_2)$ which is smooth around the point $(c(0),\eta^s_1(0))=(x,sp)$,  we obtain 
$$d(\exp_{c(0)})_{\eta_1^s(0)}\dot\eta_2^s(0)=(1-s)d(\exp_{y})_{(1-s)\overline p(0)}\dot{\overline p}(0).$$
 Moreover, 
$$\dot \eta_2^s(0)=\left[{\rm Hess}\frac{d_{CC}^2(\gamma(s),\cdot)}{2}\big|_{x}-s{\rm Hess}\frac{d_{CC}^2(y,\cdot)}{2}\big|_{x}\right] \dot c(0).$$
Finally, we recall by (\ref{wtildeY}) that $w=\dot c(0)=\overline  Y(1)\dot{\overline p}(0)$ and due to (\ref{Jacobian-Juillet}),  $\overline Y(1)$ is  invertible. 
Putting together the above computations, we have 
$$Y(s)\left[{\rm Hess}\frac{d_{CC}^2(\gamma(s),\cdot)}{2}\big|_{x}-s{\rm Hess}\frac{d_{CC}^2(y,\cdot)}{2}\big|_{x}\right]w=(1-s)\overline Y(1-s)\overline Y(1)^{-1}w.$$
Due to  
 the arbitrariness of $w$, the claim (\ref{Y-tilde-bevezetese}) follows. \hfill $\square$

\subsection{Optimal mass transportation on corank 1 Carnot groups}\label{section-omt}
We first recall some facts from Figalli and Rifford \cite{FR}. A function $\varphi:G\to \mathbb R$ is $c=d_{CC}^2/2-${\it concave} if there exist a nonempty set $S\subset G$ and a function $\varphi^c:S\to \mathbb R\cup\{-\infty\}$ with $\varphi^c \not\equiv-\infty$ such that $$\varphi(x)=\inf_{y\in S}\left\{\frac{1}{2}{d_{CC}^2(x,y)}-\varphi^c(y)\right\}.$$
If $\varphi$ is a $d_{CC}^2/2-${concave} function, let  $$\partial^c\varphi(x)=\left\{y\in S:\varphi(x)+\varphi^c(y)=\frac{1}{2}{d_{CC}^2(x,y)}\right\}$$
be the $c${\it -superdifferential of $\varphi$ at $x$}. 
For such a function $\varphi$, let $$\mathcal M_\varphi=\{x\in G:x\notin \partial^c\varphi(x)\}\ \ {\rm and}\ \ \mathcal S_\varphi=\{x\in G:x\in \partial^c\varphi(x)\}$$
be the {\it moving} and {\it static} sets, respectively. 

Let us fix $\mu_0$ and $\mu_1$  two compactly supported probability measures on $G$ which are  absolutely continuous  w.r.t. $\mathcal L^{k+1}$. According to \cite[Theorem 2.3]{FR}, there are two $d_{CC}^2/2$-concave, continuous functions $\varphi,\varphi^c:G\to \mathbb R$ such that 
\begin{equation}\label{cconcav}
\displaystyle\varphi(x)=\min_{y\in {\rm supp}(\mu_1)}\left\{\frac{1}{2}{d_{CC}^2(x,y)}-\varphi^c(y)\right\} \ {\rm and}\ 
\displaystyle\varphi^c(y)=\min_{x\in {\rm supp}(\mu_0)}\left\{\frac{1}{2}{d_{CC}^2(x,y)}-\varphi(x)\right\}
\end{equation}
and the optimal transport map is concentrated on the $c$-superdifferential of $\varphi$.
Since the distribution $\Delta$ on the corank 1 Carnot group $G$ is two-generating, it follows that $d_{CC}^2$ is locally Lipschitz on $G\times G,$  {see Agrachev and Lee  \cite[Corollary 6.2]{AL}} and Figalli and Rifford \cite[Proposition 4.2, p. 136]{FR}. Therefore, applying the version of \cite[Theorem 3.2, p. 130]{FR} with the weaker assumption for $d_{CC}^2$ of being locally Lipschitz, there exists a $d_{CC}^2/2$-concave function  $\varphi:G\to \mathbb R$ given by (\ref{cconcav}) such that $\mathcal M_\varphi$ is open and $\varphi$ is locally Lipschitz in a neighborhood of $\mathcal M_\varphi\cap {\rm supp}(\mu_0)$, thus  $\mu_0$-a.e. differentiable in $\mathcal M_\varphi$. Furthermore,  for $\mu_0$-a.e. $x$,
there exists a unique optimal transport map defined $\mu_0$-a.e. by 
{\begin{eqnarray} \label{DefIntMapH}
	\psi(x):=
	\left\{
	\begin{array}{lll}
	 \exp_x(-{\nabla}\varphi(x)) 
	\  &\mbox{if} &  x\in \mathcal M_\varphi\cap {\rm supp}(\mu_0); \\
	x &\mbox{if} &  x\in \mathcal S_\varphi\cap {\rm supp}(\mu_0),
	\end{array}\right.
	\end{eqnarray}
and for $\mu_0$-a.e. $x$ there exists a unique minimizing geodesic joining $x$ and $\psi(x)$ (or, equivalently, joining the  element $e$ with $x^{-1}\circ\psi(x)$). {Hereafter,  $\nabla=({X^0,X^{1},...,X^{d}}, Z)$ is the Carnot gradient and $\exp_x(\cdot)=x\circ \exp_e(\cdot).$}

We notice that one cannot apply directly \cite[Theorem 3.5, p. 132]{FR} of Figalli and Rifford to deduce the absolute continuity of the Wasserstein geodesic between $\mu_0$ and $\mu_1$ since in our case the semiconcavity assumption does not hold; however, we can recall the first part of their proof to conclude (based on \cite[Theorem 3.2, p. 130]{FR} and \cite[Corollary 7.22]{Villani}) that there is a unique Wasserstein geodesic $(\mu_s)_{s\in [0,1]}$ joining $\mu_0$ and $\mu_1$ given by the push-forward measure $\mu_s=(\psi_s)_\#\mu_0$ for $s\in [0,1],$ where 
{\begin{eqnarray} \label{DefIntMapH-2}
	 \psi_s(x):=
	\left\{
	\begin{array}{lll}
	 \exp_x(-s{\nabla}\varphi(x)) 
	\  &\mbox{if} &  x\in \mathcal M_\varphi\cap {\rm supp}(\mu_0); \\
	x &\mbox{if} &  x\in \mathcal S_\varphi\cap {\rm supp}(\mu_0).
	\end{array}\right.
	\end{eqnarray}
{The absolute continuity of the Wasserstein geodesic $\mu_s$ follows by the main result of Cavalletti and Mondino \cite{CM} (valid for essentially non-branching metric measure spaces) which can be state as follows:}



\begin{proposition}\label{interpolant-absolute-cont} {\rm {(Cavalletti and Mondino \cite{CM})}} Let $s\in (0,1)$.  Consider the notations introduced above and the assumptions of Theorem {\rm \ref{TJacobianDetIneq}}. Under these conditions the interpolant measure $\mu_s=(\psi_s)_{\#}\mu_0$ is  absolutely continuous w.r.t. $\mathcal L^{k+1}$.
\end{proposition}


Before  the proof of our main theorem in the next section let us indicate a 
technical difficulty that we need to address in the proof. This consists of the fact that in our setting the potential $\varphi$ generating the optimal transportation map $\psi$ via  \eqref{DefIntMapH} is not locally semiconcave {(see  Cannarsa and Sinestrari \cite{Cann-Sin})} but  only locally Lipschitz. Due to the lack of semiconcavity we do not have an Aleksandrov-type second order differentiability for 
$\varphi$ and consequently, thus we do not know  if  $\psi$ is differentiable almost everywhere. 
This regularity issue appears when we consider the transport of the mass along abnormal geodesics.

\section{Proof of the Jacobian Determinant inequality (Theorem \ref{TJacobianDetIneq})}\label{SecProof} Let $s\in (0,1).$ We shall keep the previous notations. The proof is divided into two main parts: the static and moving cases, respectively. The latter case is also divided into two parts depending how the mass is transported, i.e., along abnormal or strictly normal geodesics.  

\subsection{Static case}\label{static-case} We assume  the static set $\mathcal S_\varphi\cap {\rm supp}(\mu_0)=\{x\in {\rm supp}(\mu_0):x=\psi(x)\}$ has a positive $\mu_0$-measure. Note that $\psi_s(x)=x$ for every $x\in \mathcal S_\varphi$. If we consider the density points of $\mathcal S_\varphi$, we have that ${\rm Jac}(\psi)(x)={\rm Jac}(\psi_s)(x)=1$ for $\mu_0$-a.e. $x\in \mathcal S_\varphi$. {(Here, again Jac$(\psi)$ and Jac$(\psi_s)$ denote the  densities of $\psi_{\sharp}\mathcal L^{k+1}$ and $\psi_{s\sharp}\mathcal L^{k+1}$ w.r.t. $\mathcal L^{k+1}$.)} Note that for $x\in \mathcal S_\varphi$, we have that  $\exp_e(\theta_x)=x^{-1}\circ \psi(x)=e,$ i.e., {$\theta_x=(p_x,p_z)=(0_{\mathbb R^{k}},p_z)$ for some $p_z\in \mathbb R$}, thus $\mathcal Ap_x=0_{\mathbb R^k}$. Therefore, by the definition of the distortion coefficient, we have  $\tau_s^{k,\alpha}(\theta_x)=s$ and $\tau_{1-s}^{k,\alpha}(\theta_x)=1-s$, which concludes the proof of (\ref {Jacobi-inequality-elso}). 

\subsection{Moving case} We now assume that the moving set $\mathcal M_\varphi\cap {\rm supp}(\mu_0)$ has a positive $\mu_0$-measure.
Due to (\ref{DefIntMapH}), there exists a null $\mathcal L^{k+1}$-measure set $C_0\subset \mathcal M_\varphi\cap{\rm supp}(\mu_0)$ such that for every $x\in S:=\mathcal M_\varphi\cap {\rm supp}(\mu_0)\setminus C_0$ the function $\varphi$ is differentiable at $x$, the points $x$ and $\psi(x)$ can be joined by a unique minimizing geodesic  and $x^{-1}\circ \psi(x)=\exp_e(-\nabla\varphi(x)),$ where $$\nabla\varphi(x)=(p_x,p_z),$$
with
\begin{equation}\label{p-ix}
p_x=(X^0\varphi(x),X^1\varphi(x),...,X^d\varphi(x))\ {\rm and}\ p_z=Z\varphi(x).
\end{equation}

 Let 
\begin{equation}\label{S0set}
S_0=\{x\in S:\mathcal Ap_x= 0_{\mathbb R^{k}},\ {\rm where}\ p_x \ {\rm is\ from} \ (\ref{p-ix})\},
\end{equation}
and $$ S_1=S\setminus S_0=\{x\in S:\mathcal Ap_x\neq 0_{\mathbb R^{k}},\ {\rm where}\ p_x \ {\rm is\ from} \ (\ref{p-ix})\}.$$

We distinguish two cases.  

\subsubsection{{Moving along abnormal geodesics}}\label{moving-case-1}
 

We assume that $\mu_0(S_0)>0$. 
In terms of vector fields, the fact that $\mathcal A p_x= 0_{\mathbb R^{k}}$ with $p_x=(X^0\varphi(x),X^1\varphi(x),...,X^d\varphi(x))$ implies that $X^1\varphi(x)=...=X^d\varphi(x)=0_{\mathbb R^{2}}$ for a.e. $x\in S_0$. 
According to the explicit form of geodesics, see 
(\ref{GeodesicEqNullpz}), we have 
\begin{eqnarray}\label{abnormal-geodesic}
 \nonumber \psi_s(x)&=&x\circ \exp_e(-s\nabla \varphi(x))=x\circ \exp_e(-sX^0\varphi(x),0_{\mathbb R^{2d}},Z \varphi(x))\\  \nonumber  &=&x\circ (-sX^0\varphi(x),0_{\mathbb R^{2d+1}})\\&=& (x_1-s\partial _{x_1}\varphi(x),...,x_{k-2d}-s\partial_{x_{k-2d}} \varphi(x),x_{k-2d+1},...,x_k,z), 
\end{eqnarray} 
for a.e. $x=(x_1,...,x_k,z)\in S_0.$ In a similar way, one has 
\begin{eqnarray}\label{abnormal-geodesic-1}
 \psi(x)&=& (x_1-\partial _{x_1}\varphi(x),...,x_{k-2d}-\partial_{x_{k-2d}} \varphi(x),x_{k-2d+1},...,x_k,z), 
\end{eqnarray} 
for a.e. $x=(x_1,...,x_k,z)\in S_0.$

We divide the proof into three steps.

{{\underline{Step 1}}: \it $\varphi(\cdot, \eta, z)$ is  $d_{\mathbb R^{k-2d}}^2/2$-concave on $\mathbb R^{k-2d}$ for every $( \eta, z)\in \mathbb R^{2d}\times \mathbb R$ fixed, {i.e., for some set $S\subset \mathbb R^{k-2d}$ and function $\phi_{\eta,z}:\mathbb R^{k-2d} \to \mathbb R,$ one has $$\varphi(\xi, \eta, z)=\inf_{\overline \xi\in S}\left\{\frac{1}{2}d^2_{\mathbb R^{k-2d}}(\xi,\overline \xi)- \phi_{\eta,z}(\overline \xi)\right\}.$$}}
 Since $\varphi$ is  $d_{CC}^2/2$-concave on $G$, 
one has by (\ref{cconcav}) that for every $( \xi, \eta, z)\in \mathbb R^{k-2d}\times \mathbb R^{2d}\times \mathbb R,$
$$\varphi( \xi, \eta, z)=\min_{(\overline \xi,\overline \eta,\overline z)\in {\rm supp}(\mu_1)}\left\{\frac{1}{2}{d_{CC}^2((\xi,\eta,z), (\overline \xi,\overline \eta,\overline z))}-\varphi^c(\overline \xi,\overline \eta,\overline z)\right\}.$$

\noindent Let $\pi_1:\mathbb R^{k-2d}\times \mathbb R^{2d}\times \mathbb R\to \mathbb R^{k-2d}$ be the  projection $\pi_1( \overline \xi,\overline \eta,\overline z)=\overline \xi.$ For every $\overline \xi\in \pi_1({\rm supp}(\mu_1))$, let us introduce the compact set 
$$\Pi_{\overline \xi}=\{(\overline \eta,\overline z)\in \mathbb R^{2d}\times \mathbb R:(\overline \xi,\overline \eta,\overline z)\in {\rm supp}(\mu_1)\}.
$$ 

Let us fix $(\eta,z)\in \mathbb R^{2d}\times \mathbb R$. We notice that the function $ \phi_{\eta,z}:\pi_1({\rm supp}(\mu_1))\to \mathbb R\cup\{-\infty\}$ defined by 
$$ \phi_{\eta,z}(\overline \xi)=\max_{(\overline \eta,\overline z)\in \Pi_{\overline \xi}}\left\{\varphi^c(\overline \xi,\overline \eta,\overline z)-\frac{1}{2} {\tilde d_{CC}^2((\eta,z), (\overline \eta,\overline z))}\right\}$$
is well defined and  $\ \phi_{\eta,z}\not\equiv -\infty$. 
Since $$\displaystyle{\rm supp}(\mu_1)=\displaystyle\bigcup_{\overline \xi\in \pi_1({\rm supp}(\mu_1))}(\overline \xi,\Pi_{\overline \xi}),$$ by the Pythagorean rule (see Lemma \ref{LemmaPythagorean}) we have that for every $\xi\in \mathbb R^{k-2d}$, 
\begin{eqnarray*}
\varphi( \xi, \eta,z)&=&\min_{\overline \xi\in \pi_1({\rm supp}(\mu_1))}\min_{(\overline \eta,\overline z)\in \Pi_{\overline \xi}}\left\{\frac{1}{2}{d_{CC}^2((\xi,\eta,z), (\overline \xi,\overline \eta,\overline z))}-\varphi^c(\overline \xi,\overline \eta,\overline z)\right\}\\&=&\min_{\overline \xi\in \pi_1({\rm supp}(\mu_1))}\min_{(\overline \eta,\overline z)\in \Pi_{\overline \xi}}\left\{\frac{1}{2}d^2_{\mathbb R^{k-2d}}(\xi,\overline \xi)+\frac{1}{2} {\tilde d_{CC}^2((\eta,z), (\overline \eta,\overline z))}-\varphi^c(\overline \xi,\overline \eta,\overline z)\right\}\\&=&\min_{\overline \xi\in \pi_1({\rm supp}(\mu_1))}\left\{\frac{1}{2}d^2_{\mathbb R^{k-2d}}(\xi,\overline \xi)+\min_{(\overline \eta,\overline z)\in \Pi_{\overline \xi}}\left\{\frac{1}{2} {\tilde d_{CC}^2((\eta,z), (\overline \eta,\overline z))}-\varphi^c(\overline \xi,\overline \eta,\overline z)\right\}\right\}\\&=&\min_{\overline \xi\in \pi_1({\rm supp}(\mu_1))}\left\{\frac{1}{2}d^2_{\mathbb R^{k-2d}}(\xi,\overline \xi)- \phi_{\eta,z}(\overline \xi)\right\},
\end{eqnarray*}
which concludes the claim. 

{{\underline{Step 2}}: \it For a.e. $x=(\xi,\eta,z)\in S_0$ one can identify the Jacobian determinants ${\rm Jac}(\psi_s)(x)$ and ${\rm Jac}(\psi)(x)$ with ${\rm det}[I_{k-2d}-s{\rm Hess}_\xi(\varphi)(x)]$ and ${\rm det}[I_{k-2d}-{\rm Hess}_\xi(\varphi)(x)],$ respectively, 
	where  $I_{k-2d}$ is the $(k-2d)\times (k-2d)$ unit matrix and ${\rm Hess}_\xi(\varphi)( \xi, \eta, z)$ is the usual Euclidean Hessian of the function $\varphi(\cdot, \eta, z)$ at the point $ \xi$. }

\noindent By Step 1,  the $d_{\mathbb R^{k-2d}}^2/2$-concavity of $\varphi(\cdot, \eta, z)$  is equivalent to the convexity of $\xi\mapsto \frac{\|\xi\|_{\mathbb R^{k-2d}}^2}{2}-\varphi(\xi, \eta, z)$ on $\mathbb R^{k-2d}$. In particular, by the Aleksandrov's second differentiability theorem, the latter function is twice differentiable a.e., and its Hessian  $I_{k-2d}-{\rm Hess}_\xi(\varphi)( \xi, \eta, z)$ is positive semidefinite and symmetric for a.e. $ \xi\in \mathbb R^{k-2d}$; {the same is true for $I_{k-2d}-s{\rm Hess}_\xi(\varphi)( \xi, \eta, z)$, the latter being the convex combination of the positive semidefinite and symmetric matrices $I_{k-2d}$ and $I_{k-2d}-{\rm Hess}_\xi(\varphi)( \xi, \eta, z)$, respectively.} 
 
 By (\ref{abnormal-geodesic}) -- if it exists--  the formal Jacobian  of $\psi_s$ for a.e. $x=( \xi, \eta, z)\in S_0$ has the structure
$$\left[\begin{matrix}
A_s(x) & B_s(x)\\
0& I_{2d+1}
\end{matrix}\right],$$
where $A_s(x)=I_{k-2d}-s{\rm Hess}_\xi(\varphi)( \xi, \eta, z)$. Note however that $B_s(x)$ might not exist since we have no information on the differentiability of $\partial_i\varphi( \xi, \cdot, \cdot)$, $i\in \{1,...,k-2d\}.$ We shall explain below that the existence of $B_s(x)$ is not relevant as far as existence of the global Jacobi determinant {as density} is concerned. 

 
Observe first that due to Proposition \ref{interpolant-absolute-cont}, the interpolant measure  $\mu_s=(\psi_s)_{\#}\mu_0$ is absolutely continuous w.r.t. $\mathcal{L}^{k+1}$; let $\rho_s$ be its density function. Since the corank 1 Carnot group $(G,d_{CC},\mathcal L^{k+1})$ is a nonbranching metric measure space, both $\psi$ and $\psi_s$ are injective maps on a set of full measure of ${\rm supp}(\mu_0)$. Thus, the push-forward  measures  $\mu_s=(\psi_s)_{\#}\mu_0$ and $\mu_1=\psi_{\#}\mu_0$ and standard changes of variable should provide the  Monge-Amp\`ere equations 
\begin{equation}\label{MA-0}
\rho_0(x)=\rho_s(\psi_s(x)){\rm Jac}(\psi_s)(x) \mbox{ and } \rho_0(x)=\rho_1(\psi(x)){\rm Jac}(\psi)(x) \mbox{ for } \mu_0\mbox{-a.e. } x \in S_0.
\end{equation}
However, as we pointed out, the {differentials $D\psi$ and $D\psi_s$} may not exist on a set $S\subset S_0$ of positive measure, which requires a reinterpretation of the Monge-Amp\`ere equations in (\ref{MA-0}); we shall consider only the first term since the other one works similarly.

First of all, $\mu_s=(\psi_s)_{\#}\mu_0$  implies    
\begin{equation}\label{transport-definition}
\int_G h(y){\rm d}\mu_s(y)=\int_G h(\psi_s(x)){\rm d}\mu_0(x)
\end{equation}
for every Borel  function $h:G\to [0,\infty)$.   In particular, for every measurable set $\tilde S\subset S_0$ with positive measure and Borel function $h$ with ${\rm supp}(h) \subseteq \psi_s(\tilde S)$  we have 
$$\int_G h(y){\rm d}\mu_s(y)=\int_G h(y)\rho_s(y){\rm d}\mathcal L^{k+1}(y)=\int_{\psi_s(\tilde S)} h(y)\rho_s(y){\rm d}\mathcal L^{k+1}(y).$$
Let  $\pi_2:\mathbb R^{k-2d}\times \mathbb R^{2d}\times \mathbb R\to \mathbb R^{2d+1}$ be the projection $\pi_2(  y)=\pi_2(y_1, y_2, y_3)=(y_2, y_3)$  and for every $(y_2,y_3)\in \pi_2(\psi_s(\tilde S))$, let  $\Pi_{{(y_2,y_3)}}=\{y_1\in \mathbb R^{k-2d}:(y_1,y_2,y_3)\in \psi_s(\tilde S)\}.$  
It is clear that $\psi_s(\tilde S)=\cup_{(y_2,y_3)\in \pi_2(\psi_s(\tilde S))}(\Pi_{{(y_2,y_3)}},y_2,y_3)$; by 
Fubini's theorem it follows that 
$$\int_{\psi_s(\tilde S)} h(y)\rho_s(y){\rm d}\mathcal L^{k+1}(y)=\int_{\pi_2(\psi_s(\tilde S))} \left(\int_{\Pi_{{(y_2,y_3)}}} h(y)\rho_s(y){\rm d}\mathcal L^{k-2d}(y_1)\right){\rm d}\mathcal L^{2d+1}(y_2,y_3).$$
We consider the change of variables $y=(y_1,y_2,y_3)=\psi_s(x)$ with $x=(\xi,\eta,z)$ which shows through (\ref{abnormal-geodesic}) that $y_1=\pi_1(\psi_s(x))$ and $(y_2,y_3)=(\eta,z)$. Thus, ${\rm d}\mathcal L^{k-2d}(y_1)={\rm det}[A_s(x)]{\rm d}\mathcal L^{k-2d}(\xi)$ and $\Pi_{{(y_2,y_3)}}=\pi_1(\psi_s(\tilde S_{\eta,z},\eta,z))$ where $\tilde S_{\eta,z}=\{\xi\in \mathbb R^{k-2d}:(\xi,\eta,z)\in \tilde S\}$. Moreover, since $\pi_2(\psi_s(\tilde S))=\pi_2(\tilde S)$, the latter term  in the above relation becomes 
$$\int_{\pi_2(\tilde S)} \left(\int_{\tilde S_{\eta,z}} h(\psi_s(x))\rho_s(\psi_s(x)){\rm det}[A_s(x)]{\rm d}\mathcal L^{k-2d}(\xi)\right){\rm d}\mathcal L^{2d+1}(\eta,z)$$ which is nothing but $$\int_{\tilde S}h(\psi_s(x))\rho_s(\psi_s(x)){\rm det}[A_s(x)]{\rm d}\mathcal L^{k+1}(x).$$
The latter expression, relation  $$\int_G h(\psi_s(x)){\rm d}\mu_0(x)=\int_{\tilde S} h(\psi_s(x))\rho_0(x){\rm d}\mathcal L^{k+1}(x)$$ and (\ref{transport-definition}) together with  the arbitrariness of $h$ and $\tilde S\subset S_0$ 
 give that $$\rho_0(x)=\rho_s(\psi_s(x)){\rm det}[A_s(x)]\ \mbox{ for } \mu_0\mbox{-a.e. } x \in S_0.$$ Consequently, (\ref{MA-0}) enables us to identify $${\rm Jac}(\psi_s)(x):={\rm det}[A_s(x)]={\rm det}[I_{k-2d}-s{\rm Hess}_\xi(\varphi)( \xi, \eta, z)]\ \mbox{ for } \mu_0\mbox{-a.e. } x \in S_0.$$

{{\underline{Step 3}}: \it proof of Theorem \ref{TJacobianDetIneq} concluded $($abnormal mass transportation$)$.} 
Since $$I_{k-2d}-s{\rm Hess}_\xi(\varphi)(x)=(1-s)I_{k-2d}+s(I_{k-2d}-{\rm Hess}_\xi(\varphi)(x)),$$ we may apply the concavity of  det$(\cdot)^\frac{1}{k-2d}$ on the cone of $(k-2d)\times (k-2d)$ positive semidefinite symmetric matrices, obtaining through Step 2 that
\begin{equation}\label{Jacobian-euclidean-like}
\left[{\rm Jac}(\psi_s)(x)\right]^\frac{1}{k-2d}\geq 1-s+s\left[{\rm Jac}(\psi)(x)\right]^\frac{1}{k-2d}\ {\rm a.e.}\ x\in S_0.
\end{equation}
Now, the concavity of the function $t\mapsto t^\frac{k-2d}{k+1}, t>0,$ gives that
\begin{eqnarray*}
	\left[{\rm Jac}(\psi_s)(x)\right]^\frac{1}{k+1}&=&\left(\left[{\rm Jac}(\psi_s)(x)\right]^\frac{1}{k-2d}\right)^\frac{k-2d}{k+1}\\&\geq & \left(1-s+s\left[{\rm Jac}(\psi)(x)\right]^\frac{1}{k-2d}\right)^\frac{k-2d}{k+1}\\&\geq& 1-s+s\left[{\rm Jac}(\psi)(x)\right]^\frac{1}{k+1}\ {\rm for \ a.e.}\ x\in S_0, 
\end{eqnarray*}
which is exactly the required inequality (\ref{Jacobi-inequality-elso}).

\subsubsection{{Moving along strictly normal geodesics}}\label{moving-case-2}
We assume that $\mu_0(S_1)>0$. The proof will be divided into four steps. 

{\underline{Step 1}:} {\it $\varphi$ admits a Hessian for a.e. $x \in S_1$.} 

\noindent It is well known that the Euclidean squared distance function $d^2_{\mathbb{R}^{k-2d}}$ is semiconcave on $\mathbb{R}^{k-2d}\times \mathbb{R}^{k-2d}$, see Cannarsa and Sinestrari \cite{Cann-Sin}. Moreover, since the distribution $\tilde \Delta=\{X^1,...,X^d\}=\{X_{k-2d+1},..., X_{k}\}$ is fat on $\mathbb R^{2d+1}$, according to  Figalli and Rifford \cite[Proposition 4.1, pg. 136]{FR}, the squared distance function $\tilde{d}_{CC}^2$ is locally semiconcave on $\mathbb{R}^{2d+1} \times \mathbb{R}^{2d+1} \setminus \tilde{\mathcal{D}}$, where $\tilde{\mathcal{D}}$ denotes the diagonal of the set $\mathbb{R}^{2d+1} \times \mathbb{R}^{2d+1}$, namely $\tilde{\mathcal{D}} = \{((\eta,z),(\eta,z)) : (\eta,z) \in \mathbb{R}^{2d}\times \mathbb{R}\}$. Consequently, by the Pythagorean rule (see Lemma \ref{LemmaPythagorean}),  the squared distance function  $d^2_{CC}$ is locally semiconcave on $G \times G \setminus \mathcal{D}$, where 
\begin{equation}\label{diameter-duzzasztott}
\mathcal{D} = \{((\xi, \eta, z),(\xi', \eta, z)) : \xi, \xi' \in \mathbb{R}^{k-2d}, (\eta,z) \in \mathbb{R}^{2d}\times \mathbb{R}\}.
\end{equation}
 In order to conclude the claim, we slightly modify the  proof of \cite[Theorem 3.2]{FR}. Namely, if $x\in S_1$ is arbitrarily fixed, we have that $\mathcal Ap_x\neq 0_{\mathbb R^k}$ with $p_x$ from (\ref{p-ix}), i.e., $(x,\psi(x))\notin \mathcal D.$ In particular, if $x=(\xi_x, \eta_x, z_x)$ and $
\psi(x)=(\xi'_x, \eta'_x, z'_x)$ then  $\tilde d_{CC}((\eta_x, z_x),( \eta'_x, z'_x))=:r_x>0.$ {Due to the closeness of $\partial^c\varphi$ on $G\times G$,} there exists an open neighborhood $V_x\subset \mathcal M_\varphi\cap {\rm supp}(\mu_0)$ of $x$ such that $\tilde d_{CC}((\eta_w, z_w),( \eta'_w, z'_w))>\frac{r_x}{2}$ for every $w=(\xi_w, \eta_w, z_w)\in V_x$ and ${\psi(w) =} (\xi_w', \eta_w', z_w')\in \partial^c \varphi(w).$ Let $\tilde \varphi_{x}:G\to \mathbb R$ be defined by 
$$\tilde \varphi_{x}(w)=\inf\left\{\frac{1}{2}d_{CC}^2(w,y)-\varphi^c(y):y=(\xi_y, \eta_y, z_y)\in {\rm supp}(\mu_1),\ \tilde d_{CC}((\eta_w, z_w),( \eta_y, z_y))>\frac{r_x}{2}\right\},$$ 
where 
$w=(\xi_w, \eta_w, z_w)$. Now, the locally semiconcavity of $d_{CC}^2$ on  $G \times G \setminus \mathcal{D}$ is inherited by the $d_{CC}^2/2$-concave function $\tilde \varphi_x$  on $V_x$. Moreover, one can observe that $\tilde \varphi_x=\varphi$ on $V_x$, thus $\varphi$ is semiconcave on $V_x$.  By the Aleksandrov-Bangert theorem, see \cite[Theorem 3.10, pg. 238]{CMS},  we conclude that $\varphi$ admits a Hessian a.e. on $V_x$, concluding the claim.   

{\underline{Step 2}:} $\psi(x)\notin {\rm cut}_G(x)$ {\it for a.e.} $x\in S_1$. {We know that $\mu_0$-a.e. $x$ there exists a unique minimizing geodesic joining $x$ and $\psi(x)$, thus $\psi(x)\notin {\rm cut}_G(x)$  for a.e. $x\in S_1$.}

{\underline{Step 3}:} {\it $\psi_s$ and $\psi$ are differentiable a.e. on $S_1;$ moreover,   for a.e.  $x\in S_1,$
\begin{equation}\label{derivalt-psi-s}
d\psi_s(x)=Y_x(s)\left[{\rm Hess}\frac{d_{CC}^2(\psi_s(x),\cdot)}{2}\big|_{x}-s{\rm Hess}\varphi(x)\right],
\end{equation}
\begin{equation}\label{derivalt-psi}
d\psi(x)=Y_x(1)\left[{\rm Hess}\frac{d_{CC}^2(\psi(x),\cdot)}{2}\big|_{x}-{\rm Hess}\varphi(x)\right],
\end{equation}
 where 
 $$Y_x(s)=d(\exp_x)_{-s\nabla \varphi(x)},\ s\in (0,1].$$}

\noindent For the first part, we recall that for every $x\in S_1,$$$\psi(x)= x \circ \exp_e (-\nabla\varphi(x))\ \ {\rm and}\ \  \psi_s(x)= x \circ \exp_e(-s\nabla\varphi(x)),$$  see  
 (\ref{DefIntMapH})  and (\ref{DefIntMapH-2}), respectively. Since 
 $\psi(x) \notin {\rm cut}_G(x)$ for a.e. $x\in S_1$ (Step 2), thus $-\nabla \varphi (x)$ belongs to the injectivity domain $D$ of $\exp_e,$ and $\varphi$ has a Hessian a.e. on $S_1$ (Step 1), it follows that $\psi$ and $\psi_s$ are differentiable at a.e. $x \in S_1$. 
 
In order to prove (\ref{derivalt-psi-s}) and (\ref{derivalt-psi}), we 
need a discrete version of Claim \ref{claim-diff}:
 \begin{claim}\label{claim-diff-discrete}
 		Let $m\in \mathbb N$,  $F:\mathbb R^{2m}\to \mathbb R^m$ be a smooth function in a neighborhood of $(x,y)\in \mathbb R^{2m}$ and $\{x_n\},\{y_n\},\{z_n\}\subset \mathbb R^m$ be three sequences satisfying the following properties: 
 		\begin{enumerate}
 			\item[(a)] $\lim_{n\to \infty}x_n= x$ and $x_n\neq x$ for every $n\in \mathbb N;$   
 			\item[(b)] $\lim_{n\to \infty}y_n= y$ and $F(x_n,y_n)=F(x,y)$ for every $n\in \mathbb N;$
 			\item[(c)] $\lim_{n\to \infty}z_n= 0_{\mathbb R^{m}}$ and $\lim_{n\to \infty}\frac{z_n}{\|x_n-x\|_{\mathbb R^{m}}}= v\in \mathbb R^m.$
 		\end{enumerate}
 		Then 
 		\begin{equation}\label{discrete-kovetkeztetes}
 		\lim_{n\to \infty}\frac{F(x_n,y_n+z_n)-F(x,y)}{\|x_n-x\|_{\mathbb R^{m}}}=D_2F(x,y)v.
 		\end{equation}
 \end{claim}
 
 The proof of the claim is left as an exercise to the interested reader. 
 
 \noindent We shall apply the above claim to prove only (\ref{derivalt-psi-s}) since the proof of  (\ref{derivalt-psi}) works in a similar way. To do this, without loss of generality, we can fix a Lebesgue density point $x\in S_1$ in the differentiability set of $\varphi$, i.e., where $\varphi$ is twice differentiable (thus both $\nabla \varphi(x)$ and Hess$\varphi(x)$ exist). 
 
 Since $x$ is a Lebesgue density point of $S_1$, we can find a  linearly independent frame $\{v_i:i=1,...,k+1\}$ at $x$, such that there exist sequences $\{x_{n,i}\}\subset \mathbb R^{k+1}\setminus {\rm cut}_G(\psi_s(x))$ in the differentiability set of $\varphi$ such that for every $i\in \{1,...,k+1\}$: 
 	\begin{equation}\label{discrete-1}
 	\lim_{n\to \infty}x_{n,i}=x,\ x_{n,i}\neq x\ {\rm  for\ every}\ n\in \mathbb N,\ {\rm  and} \ \lim_{n\to \infty}\frac{x_{n,i}-x}{\|x_{n,i}-x\|_{\mathbb R^{k+1}}}=v_i;
 	\end{equation}
 	\begin{equation}\label{discrete-2}
 	\lim_{n\to \infty}\nabla \varphi (x_{n,i})=\nabla \varphi (x);
 	\end{equation}
 	$$
 	\lim_{n\to \infty}\frac{\left[\nabla\frac{d_{CC}^2(\psi_s(x),\cdot)}{2}\big|_{x_{n,i}}-s\nabla\varphi(x_{n,i})\right]-\left[\nabla\frac{d_{CC}^2(\psi_s(x),\cdot)}{2}\big|_{x}-s\nabla\varphi(x)\right]}{\|x_{n,i}-x\|_{\mathbb R^{k+1}}}=\ \ \ \ \ \ \ \ \ \ \ \ \
 	$$
  	\begin{equation}\label{discrete-3}
  	\ \ \ \ \ \ \ \ \ \ \ \ \ \ \ \ \ \ \ \ \ =\left[{\rm Hess}\frac{d_{CC}^2(\psi_s(x),\cdot)}{2}\big|_{x}-s{\rm Hess}\varphi(x)\right]v_i;\end{equation}
  \begin{equation}\label{discrete-4}
  \lim_{n\to \infty}\frac{\psi_s (x_{n,i})-\psi_s(x)}{\|x_{n,i}-x\|_{\mathbb R^{k+1}}}=d \psi_s (x)v_i.
  \end{equation}
 Fix $i\in \{1,...,k+1\}.$  We shall apply Claim \ref{claim-diff-discrete} with the smooth function $F(w,q)=\exp_w(q)$ in a neighborhood of the point $(x,y):=(x,-s\nabla \varphi(x))$ and three sequences $x_{n,i}$, $y_{n,i}:=-\nabla\frac{d_{CC}^2(\psi_s(x),\cdot)}{2}\big|_{x_{n,i}}$ and $z_{n,i}:=-y_{n,i}
  -s\nabla \varphi(x_{n,i}).$ We clearly have that $\psi_s(x_{n,i})=F(x_{n,i},y_{n,i}+z_{n,i})$. According to Proposition \ref{prop-carnot-exp}, we have that $F(x_{n,i},y_{n,i})=\psi_s(x)=F(x,y)$ for every $n\in \mathbb N$ and $\lim_{n\to \infty}y_{n,i}=-\nabla\frac{d_{CC}^2(\psi_s(x),\cdot)}{2}\big|_{x}=-s\nabla \varphi(x)=y.$ The latter relation and (\ref{discrete-2}) give that $\lim_{n\to \infty}z_{n,i}=-y+s\nabla \varphi(x)=0_{\mathbb R^{k+1}}$, while (\ref{discrete-3}) and (\ref{discrete-1}) yield  that 
  $$\lim_{n\to \infty}\frac{z_{n,i}}{\|x_{n,i}-x\|_{\mathbb R^{k+1}}}=\left[{\rm Hess}\frac{d_{CC}^2(\psi_s(x),\cdot)}{2}\big|_{x}-s{\rm Hess}\varphi(x)\right]v_i=:v\in \mathbb R^{k+1}.$$
  Thus, (\ref{discrete-kovetkeztetes}) together with (\ref{discrete-4}) reads  as 
 $$d \psi_s (x)v_i=D_2F(x,-s\nabla \varphi(x))v=d(\exp_x)_{-s\nabla \varphi(x)}\left[{\rm Hess}\frac{d_{CC}^2(\psi_s(x),\cdot)}{2}\big|_{x}-s{\rm Hess}\varphi(x)\right]v_i.$$
 Since span$\{v_1,...,v_{k+1}\}=\mathbb R^{k+1}$, the latter relation yields (\ref{derivalt-psi-s}).


 {{\underline{Step 4}}: \it proof of Theorem \ref{TJacobianDetIneq} concluded $($strictly normal mass transportation$)$.} 
%
%
We recall by Proposition \ref{prop-hessian} that the Hessian 
$$H_{x,\psi(x)}(s):={\rm Hess}\frac{d_{CC}^2(\psi_s(x),\cdot)}{2}\big|_{x}-s{\rm Hess}\frac{d_{CC}^2(\psi(x),\cdot)}{2}\big|_{x}$$
is a $(k+1)\times (k+1)$ type positive semidefinite, symmetric matrix. Since $$\nabla\frac{d_{CC}^2(\psi(x),\cdot)}{2}\big|_{x}-\nabla \varphi(x)=0_{\mathbb R^{k+1}}\ \ {\rm for\ a.e.}\ x\in S_1,$$   a similar argument as in the first part of the proof of Proposition \ref{prop-hessian} and the $d_{CC}^2/2$-concavity of $\varphi$ gives that ${\rm Hess}\frac{d_{CC}^2(\psi(x),\cdot)}{2}\big|_{x}-{\rm Hess}\varphi(x)$ is also a positive semidefinite, symmetric matrix for a.e. $x\in S_1$. 
Thus, by the concavity of  det$(\cdot)^\frac{1}{k+1}$ on the set of $(k+1)\times (k+1)$ type positive semidefinite, symmetric matrices one has\\
\\
$({\rm Jac}(\psi_s)(x))^\frac{1}{k+1} =$
\begin{eqnarray*}
	&&\quad =\det\left(Y_x(s)\left[{\rm Hess}\frac{d_{CC}^2(\psi_s(x),\cdot)}{2}\big|_{x}-s{\rm Hess}\varphi(x)\right]\right)^\frac{1}{k+1}\\
	&& \quad =
	\det(Y_x(s))^\frac{1}{k+1}\det\left[(1-s)\frac{H_{x,\psi(x)}(s)}{1-s}+s\left({\rm Hess}\frac{d_{CC}^2(\psi(x),\cdot)}{2}\big|_{x}-{\rm Hess}\varphi(x)\right)\right]^\frac{1}{k+1}\\
	&&\quad \geq 
	\det(Y_x(s))^\frac{1}{k+1}\left((1-s)\det\left[\frac{H_{x,\psi(x)}(s)}{1-s}\right]^\frac{1}{k+1}+s\det\left({\rm Hess}\frac{d_{CC}^2(\psi(x),\cdot)}{2}\big|_{x}-{\rm Hess}\varphi(x)\right)^\frac{1}{k+1}\right)\\
	&& \quad  =
	(1-s)\det(\overline Y_x(1-s)\overline Y_x(1)^{-1})^\frac{1}{k+1}+s\det( Y_x(s) Y_x(1)^{-1})^\frac{1}{k+1}({\rm Jac}(\psi)(x))^\frac{1}{k+1},
\end{eqnarray*} 
where $\overline Y_x$ corresponds to $Y_x$ via (\ref{Y-tilde-bevezetese}).

On one hand, by (\ref{Jacobi-left})  we have  that $$\det(Y_x(s) Y_x(1)^{-1})=\frac{{\rm Jac}(\exp_e)(-s\nabla \varphi(x))}{{\rm Jac}(\exp_e)(-\nabla \varphi(x))},$$ thus  by (\ref{Jacobian-Juillet}),
$$s\det(Y_x(s) Y_x(1)^{-1})^\frac{1}{k+1}=\tau_s^{k,\alpha}(\theta_x),$$ where 
$\theta_x=-\nabla \varphi(x)\in D.$ On the other hand, by the definition of $\overline Y_x$ (see (\ref{tildeY})) and relation (\ref{reverse-repres}) we also have 
$$ (1-s)\det(\overline Y_x(1-s)\overline Y_x(1)^{-1})^\frac{1}{k+1}=\tau_{1-s}^{k,\alpha}(\theta_x).$$
Combining the above facts we obtain the required Jacobian inequality (\ref{Jacobi-inequality-elso}). \hfill $\square$

\begin{remark}\rm 
 (a) Step 2 is the most fastidious part of the proof  whenever
 the mass transportation is realized along abnormal geodesics, see \S \ref{moving-case-1}.  Note that reversing the roles of the metrics $d_{\mathbb R^{k-2d}}^2$ and $\tilde d_{CC}^2$, a similar argument as in Step 1 shows that $\varphi(\xi,\cdot,\cdot)$ is a $\tilde d_{CC}^2/2-$concave function on $\mathbb R^{2d}\times \mathbb R$ ($\xi\in \mathbb R^{k-2d}$ is fixed). However, since $(x,\psi(x))\in \mathcal D$ for every $x\in S_0$ (see (\ref{S0set}) and (\ref{diameter-duzzasztott})) and we only know that $\tilde d_{CC}^2$
 is locally semiconcave on $\mathbb{R}^{2d+1} \times \mathbb{R}^{2d+1} \setminus \tilde{\mathcal{D}}$, where $\tilde{\mathcal{D}}=\pi_2(\mathcal{D})$, no conclusion can be drawn in general for the locally semiconcavity of $\varphi(\xi,\cdot,\cdot)$ on $\pi_2(S_0).$ Thus, no higher regularity is known for $\varphi(\xi,\cdot,\cdot)$  which justifies the block-decomposition of the Jacobian matrix of $\psi$ in order to interpret and compute its determinant.
 
 (b) If $S_0\subset G$ is open and $\varphi$ is smooth enough on $S_0$ (say $C^1$), one can see that  $X^1\varphi(x)=...=X^d\varphi(x)=0_{\mathbb R^{2}}$ for every $x\in S_0$ (see \S \ref{moving-case-1}) implies the fact that $\varphi$ does not depend on the components $x_{k-2d+1},.., x_k, z$, i.e., the Jacobian of $\psi$ can be calculated in the usual way on $S_0$; in particular, Example \ref{example} below falls into this framework.  
\end{remark}

We conclude this section by constructing  two measures and the optimal transportation map bet\-ween them such that a positive mass is transported along abnormal geodesics while the complementary mass is transported along strictly normal geodesics, respectively.

\begin{example}\rm \label{example}
Let $G=\mathbb R^m\times \mathbb H^d$ be the $m+2d+1$ dimensional corank 1 Carnot group endowed with its natural group operation inherited by the Euclidean space $\mathbb R^m$ and Heisenberg group $\mathbb H^d;$ in our setting,  $k=m+2d$ and $\alpha_i=4$ for every $i\in \{1,...,d\}$ in (\ref{matrix-representation}). 
   Let $a\in \mathbb R^m\setminus \{0_{\mathbb R^m}\}$ and $b\in \mathbb R^{2d}\setminus \{0_{\mathbb R^{2d}}\}$ 
 and consider the potentials $\varphi_0,\varphi_1:G\to \mathbb R$ defined by 
$$\varphi_0(x_1,x_2)=\langle a,x_1\rangle_{\mathbb R^m}\  {\rm and}\ \varphi_1(x_1,x_2)=-\langle b,{z_{x_2}}\rangle_{\mathbb R^{2d}}$$ for every 
$(x_1,x_2)=(x_1,{(z_{x_2},t_{x_2})})\in \mathbb R^m\times \mathbb H^d,$
where $\langle \cdot,\cdot \rangle_{\mathbb R^l}$ denotes the usual inner product in $\mathbb R^l.$ Moreover, let $\varphi_0^c,\varphi_1^c:G\to \mathbb R$ be defined by $$\varphi_0^c(y_1,y_2)=-\frac{1}{2}\|a\|_{\mathbb R^m}^2-\langle a,y_1\rangle_{\mathbb R^m}\  {\rm and}\ \varphi_1^c(y_1,y_2)=-\frac{1}{2}\|b\|_{\mathbb R^{2d}}^2+\langle b,{z_{y_2}}\rangle_{\mathbb R^{2d}}$$
 for every $(y_1,y_2)=(y_1,{(z_{y_2},t_{y_2})})\in \mathbb R^m\times \mathbb H^d$.
If $d_{CC}$ is the Carnot-Carath\'eodory metric  on $G$, one has for every $(x_1,x_2)\in G$ that  $$\varphi_j(x_1,x_2)=\inf_{(y_1,y_2)\in \mathbb R^m\times \mathbb H^d}\left\{\frac{1}{2}d_{CC}^2((x_1,x_2),(y_1,y_2))-\varphi_j^c(y_1,y_2)\right\},\ j\in \{0,1\},$$
where we {exploit} Lemma \ref{LemmaPythagorean}, and Ambrosio and Rigot \cite[Example 5.7, p.287]{AR} in the case $j=1$. 

Accordingly, $\varphi_j$ are $d_{CC}^2/2$-concave functions on $G$, for $j\in \{0,1\}$. 
If $\varphi=\min\{\varphi_0,\varphi_1\}$ and $\varphi^c=\max\{\varphi_0^c,\varphi_1^c\}$, we claim that   
for every $(x_1,x_2)\in G$,
\begin{equation}\label{ket-fuggveny-minimum}
\varphi(x_1,x_2)=\inf_{(y_1,y_2)\in \mathbb R^m\times \mathbb H^d}\left\{\frac{1}{2}d_{CC}^2((x_1,x_2),(y_1,y_2))-\varphi^c(y_1,y_2)\right\}.
\end{equation}  
To see this, let $(x_1,x_2)\mapsto \eta(x_1,x_2)$ be the function at the right hand side of (\ref{ket-fuggveny-minimum}). First, we  have by definition that $\varphi_j^c(y_1,y_2)\leq \varphi^c(y_1,y_2)$ for every  $(y_1,y_2)\in G$ and $j\in \{0,1\}.$ Accordingly, $\varphi_j(x_1,x_2)\geq \eta(x_1,x_2)$ for every  $(x_1,x_2)\in G$ and $j\in \{0,1\},$ i.e., $\varphi\geq \eta$. 

To check the converse inequality, {we provide a generic argument, independent from the Carnot structure}. Fix $(x_1,x_2)\in G$ arbitrarily and assume {without loss of generality} that $\varphi_0(x_1,x_2)\leq \varphi_1(x_1,x_2)$. Then for every $(y_1,y_2)\in G$, we have 
$$\varphi_0(x_1,x_2)\leq \frac{1}{2}d_{CC}^2((x_1,x_2),(y_1,y_2))-\varphi_0^c(y_1,y_2);$$
$$\varphi_0(x_1,x_2)\leq \varphi_1(x_1,x_2)\leq \frac{1}{2}d_{CC}^2((x_1,x_2),(y_1,y_2))-\varphi_1^c(y_1,y_2).$$
Consequently, for every $(y_1,y_2)\in G$, one has 
\begin{eqnarray*}
\varphi_0(x_1,x_2)&\leq&\frac{1}{2}d_{CC}^2((x_1,x_2),(y_1,y_2))+{\min}\{-\varphi_0^c(y_1,y_2),-\varphi_1^c(y_1,y_2)\}\\&=&
\frac{1}{2}d_{CC}^2((x_1,x_2),(y_1,y_2))-{\max}\{\varphi_0^c(y_1,y_2),\varphi_1^c(y_1,y_2)\}\\&=&
\frac{1}{2}d_{CC}^2((x_1,x_2),(y_1,y_2))-\varphi^c(y_1,y_2).
\end{eqnarray*}
Taking the infimum on the right w.r.t. $(y_1,y_2)\in G$, we obtain 
$\varphi_0(x_1,x_2)\leq \eta(x_1,x_2),$ i.e., $\varphi(x_1,x_2)\leq \eta(x_1,x_2)$, which concludes the proof of (\ref{ket-fuggveny-minimum}). In particular, (\ref{ket-fuggveny-minimum}) implies that  
 $\varphi$ is a $d_{CC}^2/2$-concave function on $G$.
 
  Let $G^0=\{(x_1,x_2)=(x_1,(z_2,t))\in \mathbb R^m\times \mathbb H^d:\langle (a,b),(x_1,z_2)\rangle_{\mathbb R^m\times \mathbb R^{2d}}=0 \}$ be the hyperplane separating $\mathbb R^m\times \mathbb R^{2d+1}$ into two halfspaces  $G^-=\{ (x_1,(z_2,t))\in \mathbb R^m\times \mathbb H^d:\langle (a,b),(x_1,z_2)\rangle_{\mathbb R^m\times \mathbb R^{2d}}\leq 0\}$
and $G^+=G\setminus G^-$. It follows that 
$$\varphi(x_1,x_2)=\left\{
\begin{array}{lll}
\varphi_0(x_1,x_2) &\mbox{if} &  (x_1,x_2)\in G^-;\\
\varphi_1(x_1,x_2) &\mbox{if} &  (x_1,x_2)\in G^+,
\end{array}
\right.$$
and $\varphi$ is  differentiable on $G\setminus G^0.$  Let $\psi:G\to G$  be the optimal transport map generated by the $d_{CC}^2/2$-concave function $\varphi$, {see Figalli and Rifford \cite{FR}}; by Proposition \ref{proposition-geodetikus} we have for every  $(x_1,x_2)\in G\setminus G^0$ that 
 $$\psi(x_1,x_2)=\exp_{(x_1,x_2)}(-\nabla \varphi(x_1,x_2))= \left\{
 \begin{array}{lll}
 \psi_0(x_1,x_2) &\mbox{if} &  (x_1,x_2)\in G^-\setminus G^0;\\
 \psi_1(x_1,x_2) &\mbox{if} &  (x_1,x_2)\in G^+,
 \end{array}
 \right.$$
 where $$\psi_0(x_1,x_2)=\exp_{(x_1,x_2)}(-a,0_{\mathbb R^{2d+1}})=(x_1-a,x_2)$$ and $$\psi_1(x_1,x_2)=\exp_{(x_1,x_2)}(0_{\mathbb R^{m}},b,0)=(x_1,x_2* (b,0));$$ here, $'*'$  denotes the group operation  on $\mathbb H^d$. 
 
Let $\mu_0=\mathbbm{1} _{B}\mathcal L^{m+2d+1}$, where $B\subset G$ is a closed ball centered at $0_{\mathbb R^{m+2d+1}}$ with  $\mathcal L^{m+2d+1}(B)=1$ and $\mu_1=\psi_\#\mu_0$. Note that every element of $B $ belongs to the moving set $M_\varphi$. Moreover, one can see that 
${\rm supp}(\mu_1)=\overline{\psi_0(B\cap G^-\setminus G^0)}\bigcup \overline{\psi_1(B\cap G^+)},$
and the sets $S_0$ and $S_1$ appearing in the proof of Theorem \ref{TJacobianDetIneq} correspond  to the two half balls $B\cap G^-$ and $B\cap G^+$ (up to  null measure sets), respectively. Consequently, the optimal mass transportation map $\psi$ translates the mass from $S_0$ along abnormal (Euclidean) geodesics into  $\tilde S_0=(-a,0_{\mathbb R^{m+2d+1}})+S_0$, while the mass from $S_1$ is transported along strictly normal (Heisenberg) geodesics into a distorted half ball $\tilde S_1=\{(x_1,x_2* (b,0)):(x_1,x_2)\in B\cap G^+\}$, see Figure \ref{abra-elso}. 

\begin{figure}[H]
	\includegraphics[scale=0.5]
	{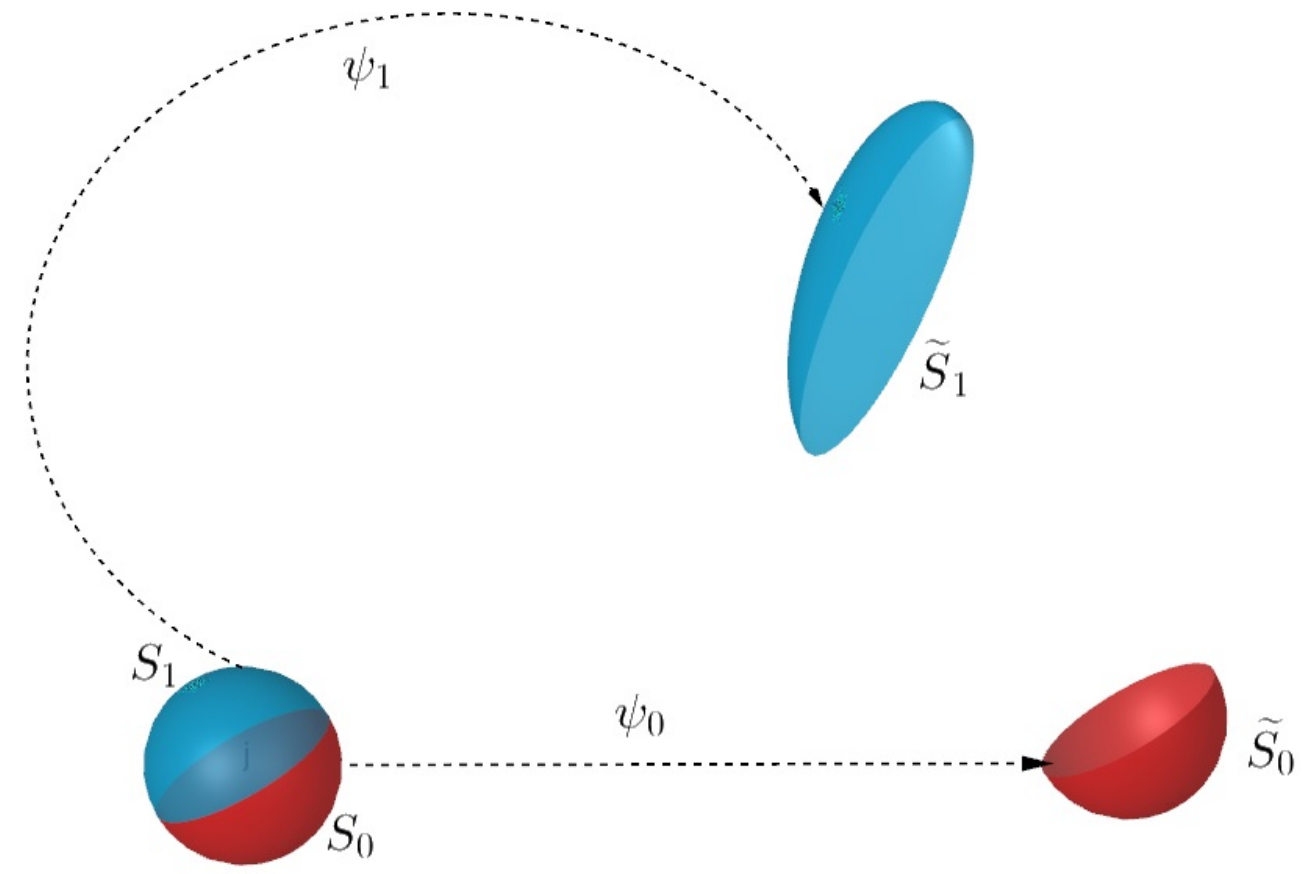}
	\caption{The half balls $S_0$ and $S_1$ are transported along abnormal and strictly normal  geodesics into the sets $\tilde S_0$ and $\tilde S_1$, respectively.}\label{abra-elso}
\end{figure}
\end{example}

\medskip

\section{Applications}\label{SecApps}

Having the Jacobian determinant  inequality (\ref{Jacobi-inequality-elso}), we can prove several functional and geometric inequalities on corank 1 Carnot groups. 


Let us denote by $\rho_0, $ $\rho_1$ and $\rho_s$ the density functions (w.r.t. $\mathcal{L}^{k+1}$) of the absolutely continuous, compactly supported measures $\mu_0$,  $\mu_1=\psi_{\#}\mu_0$ and $\mu_s=(\psi_s)_{\#}\mu_0$, respectively. In fact, we have the 
 Monge-Amp\`ere equations
\begin{equation}\label{MA}
\rho_0(x)=\rho_s(\psi_s(x)){\rm Jac}(\psi_s)(x) \mbox{ and } \rho_0(x)=\rho_1(\psi(x)){\rm Jac}(\psi)(x) \mbox{ for } \mu_0\mbox{-a.e. } x \in G.
\end{equation}
These equations can be deduced in a standard way both in the static case (see \S \ref{static-case}) and moving case with optimal mass transport along strictly normal geodesics (see \S \ref{moving-case-2}), while in the case of abnormal transportation we provided a proper interpretation of them (see \S \ref{moving-case-1}, Step 2).

Due to (\ref{MA}) we may reformulate the Jacobian determinant inequality (\ref{Jacobi-inequality-elso}) as
\begin{eqnarray}\label{Jacobi-inequality-masodik}
\left(\rho_s(\psi_s(x)\right)^{-\frac{1}{k+1}}\geq \tau_{1-s}^{k,\alpha}(\theta_x)\left(\rho_0(x)\right)^{-\frac{1}{k+1}}+\tau_{s}^{k,\alpha}(\theta_x)\left(\rho_1(\psi(x)\right)^{-\frac{1}{k+1}},
\end{eqnarray}
which holds $\mu_0$ a.e. on the restricted set $G_0 = \{x \in G : \rho_0(x) > 0\}$. Observe that by definition $G_0$ is of full measure in ${\rm supp}(\mu_0)$. For a fixed $s \in (0,1)$ we restrict $G_0$ to the injectivity domain of $\psi$ and $\psi_s$ which will be still of full measure in ${\rm supp}(\mu_0)$. 
Moreover, we may exclude those points $x\in  S_1$ from  $G_0$ for which $x^{-1} \circ \psi(x) \in {\rm cut}_G(e)$, see Step 2 in \S \ref{moving-case-2}, still obtaining a full measure set in ${\rm supp}(\mu_0)$ which prevents the blow-up of  coefficients $\tau_{1-s}^{k,\alpha}(\theta_x)$ and  $\tau_s^{k,\alpha}(\theta_x)$, respectively. 


\subsection{Entropy inequalities}
Let $(G,\circ)$ be a $k+1$ dimensional corank 1 Carnot group and  $U:[0, \infty) \to \mathbb R$ be a function. The $U$-entropy  of an absolutely continuous measure $\mu$ w.r.t. $\mathcal L^{k+1}$ on $G$ is
defined by $\displaystyle{\rm Ent}_U(\mu | \mathcal L^{k+1}) = \int_{G} U\left( \rho(x) \right) \dd \mathcal L^{k+1}(x),$
where $\rho=\frac{\dd \mu}{\dd \mathcal L^{k+1}}$ is the density function of $\mu.$

{By using the injectivity of $\psi_s$ and $\psi$ on $G_0$ (with a suitable change of variables), a similar  argument as in \cite{BKS} provides the following entropy inequality.}  

\begin{theorem}\label{TEntIneqCarnotGen} {\bf (Entropy inequality)}
Under the same assumptions as in Theorem \ref{TJacobianDetIneq}, if 
	 $U: [0, \infty) \to \mathbb R$ is a function such that $U(0)=0$ and $t \mapsto t^{k+1} U\left(\frac{1}{t^{k+1}}\right)$ is non-increasing and convex, the following entropy inequality holds:
	\begin{eqnarray*}
		{\rm Ent}_{U}(\mu_s | \mathcal L^{k+1})   &\leq&  (1-s) \int_{G} \left(\tilde{\tau}_{1-s}^{k,\alpha}(\theta_x)\right)^{k+1} U\left(\frac{\rho_0(x)}{\left(\tilde{\tau}_{1-s}^{k,\alpha}(\theta_x)\right)^{k+1}}\right) \dd \mathcal L^{k+1}(x) \\
		&&+ s \int_{G} \left(\tilde{\tau}_{s}^{k,\alpha}(\theta_{\psi^{-1}(y)})\right)^{k+1} U\left(\frac{\rho_1(y)}{\left(\tilde{\tau}_{s}^{k,\alpha}(\theta_{\psi^{-1}(y)})\right)^{k+1}}\right) \dd \mathcal L^{k+1}(y),
	\end{eqnarray*}
	where $\tilde{\tau}_{s}^{k,\alpha}=s^{-1}{\tau}_{s}^{k,\alpha}.$

\end{theorem}

\begin{corollary}\label{CEntIneqCarnotUnif} Under the same assumptions as in Theorem \ref{TEntIneqCarnotGen}, we have the following uniform entropy inequality: 
	\begin{eqnarray*}
		{\rm Ent}_{U}(\mu_s | \mathcal L^{k+1})   &\leq&  (1-s)^3 \int_{G}  U\left(\frac{\rho_0(x)}{(1-s)^2}\right) \dd \mathcal L^{k+1}(x) + s^3 \int_{G}  U\left(\frac{\rho_1(y)}{s^2}\right) \dd \mathcal L^{k+1}(y).
	\end{eqnarray*}
\end{corollary}

{\it Proof.} Since $t\mapsto \frac{\mathbb d_i(t,s)}{\mathbb d_i(t,1)}$ is increasing on $(0,2\pi)$ for every $s\in (0,1)$, $i\in \{1,2\},$ and $\lim_{t\to 0}\frac{\mathbb d_1(t,s)}{\mathbb d_1(t,1)}=1$, $\lim_{t\to 0}\frac{\mathbb d_2(t,s)}{\mathbb d_2(t,1)}=s^2$, see \cite[Lemma 2.1]{BKS}, we obtain 
\begin{eqnarray} \label{tau-lbound}
\tau^{k,\alpha}_s(\theta_x)\geq s^{\frac{k+3}{k+1}} \mbox{ for all } s \in (0,1) \mbox{ and } x \in G_0.
\end{eqnarray}
Thus, for the weights $\tilde{\tau}_s^{k, \alpha}$ we obtain 
\begin{eqnarray}\label{tau-tilde-lbound}
\left(\tilde{\tau}^{k,\alpha}_s(\theta_x)\right)^{k+1}\geq s^2 \mbox{ for all } s \in (0,1) \mbox{ and } x \in G_0.
\end{eqnarray}
Since the  map $t \mapsto t^{k+1} U\left(\frac{1}{t^{k+1}}\right)$ is non-increasing, the desired inequality directly follows from Theorem \ref{TEntIneqCarnotGen}.
\hfill $\square$

\begin{remark}\rm As a particular case of Theorem \ref{TEntIneqCarnotGen} and Corollary \ref{CEntIneqCarnotUnif}, we may choose various particular entropies for $U$, as the R\'enyi-type entropy,  Shannon entropy,  kinetic-type entropy or Tsallis entropy.
\end{remark}

\noindent 	


\subsection{Brunn-Minkowski inequalities} 
Let $(G,\circ)$ be a connected, simply connected nilpotent Lie group of (topological) dimension $N$, and  $\mu$ be a Haar measure on $G$. By extending a result of Leonardi and Masnou \cite{LM} from Heisenberg groups, Tao \cite{T} proved that for every nonempty and bounded open sets $A, B\subset G$ the multiplicative Brunn-Minkowski inequality holds on $(G,\circ)$: 
\begin{eqnarray}\label{multiplicative-BM} 
\mu(A \circ B)^\frac{1}{N} \geq \mu(A)^\frac{1}{N} + \mu(B)^\frac{1}{N}.
\end{eqnarray}
In particular,  this inequality is also valid on any $k+1$ dimensional corank 1 Carnot group $G$ with $N=k+1$ and $\mu=\mathcal L^{k+1}$. 

In the sequel, we prove geodesic Brunn-Minkowski inequalities on corank 1 Carnot groups. To do this, let $A,B \subset G$ be two nonempty sets. In the sequel we want to quantify the Carnot distortion coefficients which characterize the  sets $A$ and $B$.  For this reason we introduce the notations 
\begin{eqnarray}\label{disztrozio-A-B}
\tau_s^{k,\alpha}(A,B) = \sup_{A_0, B_0} \inf_{(x,y) \in A_0\times B_0} \{\tau_s^{k,\alpha}(p): \exp_e(p) = x^{-1}\circ y\}
\end{eqnarray}
and 
\begin{eqnarray}\label{disztrozio-A-B-1}
\tilde{\tau}_s^{k, \alpha}(A,B) = \sup_{A_0, B_0} \inf_{(x,y) \in A_0\times B_0} \{\tilde{\tau}_s^{k, \alpha}(p): \exp_e(p) = x^{-1}\circ y\} = s \tau_s^{k,\alpha}(A,B),
\end{eqnarray}
where $A_0$ and $B_0$ are nonempty, full measure subsets of $A$ and $B$, respectively.
{Note that by taking sets $A_0,B_0$ with the above properties we might obtain better coefficients than if simply take the initial sets $A,B$; more precisely, one always has 
$$\tilde{\tau}_s^{k, \alpha}(A,B)\geq \inf_{(x,y) \in A\times B} \{\tilde{\tau}_s^{k, \alpha}(p): \exp_e(p) = x^{-1}\circ y\},$$
with possibly strict inequality e.g. when some discrete points  $x\in A$ and $y\in B$ are in a particular position as $\exp_e(p) = x^{-1}\circ y$ with $p\in D$ and $p_z=0$.}
Recalling relation (\ref{reverse-repres}) between the parameters of the exponential map joining $e$ to $x \in G$ and $x^{-1}\in G$, respectively,  the following symmetry properties hold:
\begin{eqnarray}\label{tau-symm}
\tau_s^{k,\alpha}(x,y) = \tau_s^{k,\alpha}(y,x) \mbox{ and } \tilde{\tau}_s^{k, \alpha}(x,y) = \tilde{\tau}_s^{k, \alpha}(y,x) \mbox{ for all } x,y \in G.
\end{eqnarray}


For every $s\in [0,1]$ and  $x,y\in G,$ the set of $s$-intermediate points between $x$ and $y$ is
\begin{eqnarray}\label{DZs}
Z_s(x,y)=  \{ z \in G : d_{CC}(x,z) = s d_{CC}(x,y),\
d_{CC}(z,y) = (1-s) d_{CC}(x,y)\}.
\end{eqnarray}\label{Z-antisymm} 
We clearly have the antisymmetry property 
$$
Z_s(x,y) = Z_{1-s}(y,x) \mbox{ for all } x,y \in G \mbox{ and } s \in [0,1].
$$
The notion of $s$-intermediate points can be extended to the nonempty sets $A,B \subset G$ as
$$Z_s(A,B) = \bigcup_{(x,y) \in A \times B} Z_s(x,y).$$


\begin{theorem}\label{TBrunn-Minkowski-2}  {\bf (Weighted and non-weighted Brunn-Minkowski inequalities)} Let $(G,\circ)$ be a $k+1$ dimensional corank 1 Carnot group,  $s\in (0,1),$ and $A$ and $B$ be two nonempty measurable sets of $G$. Then the following inequalities hold:
	\begin{itemize}
		\item[{\rm (i)}] $\displaystyle \mathcal L^{k+1}(Z_s(A,B))^\frac{1}{k+1} \geq
		\tau_{1-s}^{k,\alpha}(A,B)\mathcal
		L^{k+1}(A)^\frac{1}{k+1}+\tau_s^{k,\alpha}(A,B)\mathcal
		L^{k+1}(B)^\frac{1}{k+1};$
		\item[{\rm (ii)}] $\displaystyle \mathcal L^{k+1}(Z_s(A,B))^\frac{1}{k+1} \geq
		(1-s)^\frac{k+3}{k+1}\mathcal
		L^{k+1}(A)^\frac{1}{k+1}+s^\frac{k+3}{k+1}\mathcal
		L^{k+1}(B)^\frac{1}{k+1};$
		\item[{\rm (iii)}] $\displaystyle \mathcal L^{k+1}(Z_s(A,B))^\frac{1}{k+3} \geq
		\left(\frac{1}{4}\right)^{\frac{1}{k+3}} \left( (1-s)\mathcal
		L^{k+1}(A)^\frac{1}{k+3}+s\mathcal
		L^{k+1}(B)^\frac{1}{k+3}\right).$
	\end{itemize}
\end{theorem}

{\it Proof.} First of all, we notice that if $Z_s(A,B)$ is not measurable,  $\mathcal L^{k+1}(Z_s(A,B))$ will denote the outer Lebesgue measure of $Z_s(A,B)$. 

(i) We first assume that both $A$ and $B$ have finite $\mathcal
L^{k+1}$-measures. If both sets have null measure, we have nothing to prove; thus, we may assume that $\max\left\{\mathcal
L^{k+1}(A),\mathcal
L^{k+1}(B)\right\}>0.$
The proof is divided into three steps. 

{{\underline{Step 1}}:
	{\it one has $\tau_s^{k,\alpha}(A,B)<\infty$ and $\tau_{1-s}^{k,\alpha}(A,B)<\infty.$} By  (\ref{disztrozio-A-B}), if $\tau_s^{k,\alpha}(A,B)=+\infty$, we have in particular that $x^{-1} \circ y \in {\rm cut}_G(e)$ for a.e. $(x,y)\in A\times B.$ Therefore, $0=\mathcal L^{k+1}({\rm cut}_G(e))\geq \mathcal L^{k+1}(A^{-1}\circ B)$. Thus, by the  multiplicative Brunn-Minkowski inequality (\ref{multiplicative-BM}) it follows that $\mathcal
	L^{k+1}(A)= \mathcal
	L^{k+1}(B)=0,$ which contradicts our initial assumption.

	{{\underline{Step 2}}:  {\it the case} $\mathcal
		L^{k+1}(A)\neq 0\neq \mathcal
		L^{k+1}(B).$ Let $\mu_0=\frac{\mathbbm{1} _A(x)}{\mathcal
			L^{k+1}(A)}\mathcal		L^{k+1}$,  $\mu_1=\frac{\mathbbm{1} _B(x)}{\mathcal
			L^{k+1}(B)}\mathcal L^{k+1}$ and the R\'enyi entropy $U(r)=-r^{1-\frac{1}{k+1}}$  $(r\geq 0)$ in Theorem \ref{TEntIneqCarnotGen}; thus the entropy inequality and relations (\ref{disztrozio-A-B}) and (\ref{disztrozio-A-B-1}) imply that
		\begin{eqnarray*}
		\int_{\psi_s(A)}\rho_s(z)^{1-\frac{1}{k+1}}\dd \mathcal{L}^{k+1}(z)&\geq& \tau_{1-s}^{k,\alpha}(A,B)\int_{A}\rho_0^{1-\frac{1}{k+1}}\dd \mathcal{L}^{k+1}+\tau_s^{k,\alpha}(A,B)\int_{B}\rho_1^{1-\frac{1}{k+1}}\dd \mathcal{L}^{k+1}\\&=&\tau_{1-s}^{k,\alpha}(A,B)\mathcal
		L^{k+1}(A)^\frac{1}{k+1}+\tau_s^{k,\alpha}(A,B)\mathcal
		L^{k+1}(B)^\frac{1}{k+1}.
		\end{eqnarray*}
	By H\"older's inequality one has that 
		\begin{eqnarray*}
	\int_{\psi_s(A)}\rho_s(z)^{1-\frac{1}{k+1}}\dd \mathcal{L}^{k+1}(z)&\leq& \left(\int_{\psi_s(A)}\rho_s(z)\dd \mathcal{L}^{k+1}(z)\right)^{1-\frac{1}{k+1}}\left(\int_{\psi_s(A)}\dd \mathcal{L}^{k+1}(z)\right)^{\frac{1}{k+1}}\\&=& \mathcal L^{k+1}(\psi_s(A))^\frac{1}{k+1}.
		\end{eqnarray*}
	Since $\psi_s(A)\subset Z_s(A,B)$, the claim follows. 
	 
		{{\underline{Step 3}}:  {\it the case} $\mathcal
			L^{k+1}(A)= 0\neq \mathcal
			L^{k+1}(B)$ {\it or} $\mathcal
			L^{k+1}(A)\neq 0= \mathcal
			L^{k+1}(B).$    
			%
		In fact, our claim reduces to proving that 	for every $x\in G$, we have 
			\begin{equation}\label{Belso-biz-Juillet}
			\displaystyle \mathcal L^{k+1}(Z_s(\{x\},B)) \geq
			\left(\tau_s^{k,\alpha}(\{x\},B)\right)^{k+1}\mathcal
			L^{k+1}(B).
			\end{equation}
			The latter inequality follows by an approximation argument. In fact, if $\{\eps_n\}_{n\in \mathbb N}$ is a decreasing sequence converging to 0, by Step 2  we have for every $n\in \mathbb N$ that
			$$\displaystyle \mathcal L^{k+1}(Z_s(B(x,\eps_{n}),B))^\frac{1}{k+1} \geq
			\tau_{1-s}^{k,\alpha}(B(x,\eps_n),B)
			\eps_n^\frac{k+2}{k+1}+
			\tau_{s}^{k, \alpha}(B(x,\eps_n),B)
			\mathcal L^{k+1}(B)^\frac{1}{k+1},$$
			where $B(x,r)=\{y\in G:d_{CC}(x,y)\leq r\}.$
			By using the monotone convergence theorem one can prove that 
		$$	\lim_{n\to \infty}\mathcal L^{k+1}(Z_s(B(x,\eps_{n}),B))=\mathcal L^{k+1}(Z_s(\{x\},B))\ \ {\rm and}\ \  	\lim_{n\to \infty}\tau_{s}^{k, \alpha}(B(x,\eps_n),B)=\tau_{s}^{k, \alpha}(\{x\},B),$$
		which proves (\ref{Belso-biz-Juillet}). 
If $A$ or $B$ has infinite $\mathcal L^{k+1}$-measure, we apply again an approximation argument.

			(ii) This property follows by (i) combined with the universal lower bound (\ref{tau-lbound}) for $\tau_s^{k,\alpha}$.

		{	(iii) Property (ii) is combined with the $p$-mean inequality (\ref{MspIneq}) below 	with the choices $a = (1-s)^{-2} $, $b = s^{-2} $, $p = \frac{1}{2}$, $q = \frac{1}{k+1}$ and $\eta = \frac{1}{k+3}$, respectively.} \hfill $\square$\\
			
			The main result of Rizzi \cite{Rizzi} concerning the measure contraction property on corank 1 Carnot groups is a direct consequence of the Brunn-Minkowski inequality (Theorem \ref{TBrunn-Minkowski-2}):

		\begin{corollary}\label{TMCP-1}  {\bf (Measure contraction property)}
			Let $(G,\circ)$ be a $k+1$ dimensional corank 1 Carnot group.  Then the  measure contraction property {\rm{\textsf{ MCP}}}$(0,k+3)$ holds
			on $G$, i.e., for every  $s\in [0,1]$, $x\in G$
			and nonempty measurable set $E\subset G$,
			$$\displaystyle \mathcal L^{k+1}(Z_s(\{x\},E)) \geq \left(\tau_s^{k,\alpha}(\{x\},E)\right)^{k+1} \mathcal L^{k+1}(E) \geq s^{k+3}\mathcal L^{k+1}(E).$$
		\end{corollary}
		\vspace{0.5cm}


\subsection{Borell-Brascamp-Lieb inequalities}

In order to formulate our Borell-Brascamp-Lieb inequalities we introduce the notion of the $p$-mean, which for two non-negative numbers $a,b$ and weight $s \in (0,1)$ is defined as
$$M_s^p(a,b)=\left\{\begin{array}{lll}
	\left( (1-s)a^p + s b^p \right)^{1/p} &\mbox{if} &  ab\neq 0, \\
	0 &\mbox{if} &  ab=0,
	\end{array}\right.$$
{with the  conventions $M_s^{-\infty}(a,b)=\min\{a,b\}$; 
$M_s^{0}(a,b)=a^{1-s}b^s;$ and $M_s^{+\infty}(a,b)=\max\{a,b\}$ if $ab\neq 0$ and $M_s^{+\infty}(a,b)=0$ if $ab= 0$.}
According to Gardner \cite[Lemma 10.1]{Gar}, one has
\begin{eqnarray}\label{MspIneq}
		M_s^{p}(a,b)M_s^{q}(c,d) \geq M_s^{\eta}(ac, bd),
\end{eqnarray} 
for every $a,b,c,d \geq 0, s \in (0,1)$ and $p, q \in \R$ such that $p+q\geq0$ with ${\eta}=\frac{pq}{p+q}$ when $p$ and $q$ are not both zero, and $\eta=0$ if $p=q=0$. 

{Having the Jacobian determinant inequality (\ref{Jacobi-inequality-masodik}), we can prove Borell-Brascamp-Lieb-type inequalities on corank 1 Carnot groups.  In the sequel we state some of them. We refer to \cite{BKS} for similar results with detailed proofs in the setting of the Heisenberg groups:} 

\begin{theorem}\label{TRescaledBBLWithWeights} {\bf (Weighted Borell-Brascamp-Lieb inequality)}
		Fix $s\in (0,1)$ and $p \geq -\frac{1}{k+1}$. 
		Let $f,g,h: G\to [0,\infty)$ be Lebesgue integrable
		functions with the property that for all $(x,y)\in G\times G, z\in Z_s(x,y),$
		\begin{eqnarray}\label{ConditionRescaledBBLWithWeights}
		h(z) \geq M^{p}_s
		\left(\frac{f(x)}{\left(\tilde{\tau}_{1-s}^{k, \alpha}(y,x)\right)^{k+1}},\frac{g(y)}{\left(\tilde{\tau}_s^{k, \alpha}(x,y)\right)^{k+1}} \right).
		\end{eqnarray}
		Then the following inequality holds:
		\begin{eqnarray*}
			\int_{G} h \geq M^\frac{p}{1+(k+1)p}_s \left(\int_{G}
			f, \int_{G} g \right).
		\end{eqnarray*}
	\end{theorem}

\begin{remark}\label{RStrongerBBL} \rm  Observe that Theorem \ref{TRescaledBBLWithWeights} holds as well under weaker conditions, namely, if inequality (\ref{ConditionRescaledBBLWithWeights}) holds only for those $x,y \in G$ for which $f(x)>0$ and $g(y)>0$.
\end{remark}

As a direct consequence of Theorem \ref{TRescaledBBLWithWeights}, inequality (\ref{tau-tilde-lbound}) and the monotonicity of the $p$-mean we can formulate the following weaker Borell-Brascamp-Lieb-type inequality:

\begin{corollary}\label{CLighterWeightedBBL}  {\bf (Uniformly weighted Borell-Brascamp-Lieb inequality)}
				Fix $s\in (0,1)$ and $p \geq -\frac{1}{k+1}.$  Let $f,g,h:
				G\to [0,\infty)$ be  Lebesgue integrable
				functions satisfying
				\begin{eqnarray}&\label{1-ConditionRescaledBBLWithoutWeights}
				h(z) \geq M^{p}_s \left(\frac{f(x)}{(1-s)^2},\frac{g(y)}{s^2}\right) \ \  for\ all\ (x,y)\in G\times
				G, z\in Z_s(x,y).
				\end{eqnarray}
				Then the following inequality holds:
				\begin{eqnarray*}
					\int_{G} h \geq M^\frac{p}{1+(k+1)p}_s \left(
					\int_{G} f, \int_{G}g\right).
				\end{eqnarray*}
			\end{corollary}
						
\begin{corollary}\label{CRescaledBBLWithoutWeights}  {\bf (Non-weighted Borell-Brascamp-Lieb inequality)}
				Fix $s\in (0,1)$ and $p \geq -\frac{1}{k+3}.$  Let $f,g,h:
				G\to [0,\infty)$ be  Lebesgue integrable
				functions satisfying
				\begin{eqnarray}&\label{ConditionRescaledBBLWithoutWeights}
				h(z) \geq M^{p}_s (f(x),g(y)) \ \  for\ all\ (x,y)\in G\times
				G, z\in Z_s(x,y).
				\end{eqnarray}
				Then the following inequality holds:
				\begin{eqnarray}\label{correction}
				\int_{G} h \geq \frac{1}{4}M^\frac{p}{1+(k+3)p}_s \left(
				\int_{G} f, \int_{G}g\right).
				\end{eqnarray}
			\end{corollary}
			
{\it Proof.}  By the $p$-mean inequality (\ref{MspIneq}) and  assumption (\ref{ConditionRescaledBBLWithoutWeights}), we have  
\begin{eqnarray}
4h(z) = M_s^p(f(x),g(y)) M_s^{\frac{1}{2}}\left(\frac{1}{(1-s)^2}, \frac{1}{s^2}\right) \geq M_s^{\frac{p}{2p+1}}\left(\frac{f(x)}{(1-s)^2}, \frac{g(y)}{s^2} \right),
\end{eqnarray}
for every $x,y \in G$ and $z \in Z_s(x,y)$.
By the assumption $p\geq -\frac{1}{k+3}$ we have $\frac{p}{2p+1}\geq-\frac{1}{k+1}$, so we can apply Corollary \ref{CLighterWeightedBBL} for the setting $\tilde{h} = 4 h$, $\tilde{f} = f$, $\tilde{g} = g$ and $\tilde{p} = \frac{p}{2p+1}$,   obtaining the desired inequality.	\hfill $\square$

\begin{remark} \rm (a)
All three versions of the Borell-Brascamp-Lieb inequality imply a corresponding Pr\'ekopa-Leindler-type inequality by setting $p=0$ and using the convention $M_s^0(a,b) = a^{1-s}b^s$ for all $a,b \geq 0$ and $s \in (0,1)$. 


(b)  The Brunn-Minkowski inequality (i) in Theorem \ref{TBrunn-Minkowski-2} can be obtained alternatively from Theorem \ref{TRescaledBBLWithWeights} whenever  $\mathcal
	L^{k+1}(A)\neq 0\neq \mathcal
	L^{k+1}(B).$ Indeed, let $p=+\infty$, and choose the functions  $f(x)=\left(\tilde{\tau}_{1-s}^{k, \alpha}(A,B)\right)^{k+1} \mathbbm{1} _A(x),$
$g(y)=\left(\tilde{\tau}_{s}^{k, \alpha}(A,B)\right)^{k+1} \mathbbm{1} _B(y)$ and $h(z)=\mathbbm{1}
_{Z_s(A,B)}(z).$ With these choices assumption 
(\ref{ConditionRescaledBBLWithWeights}) holds at the points $x,y \in G$ where $f(x)>0$ and $g(y)>0$ and due to Remark \ref{RStrongerBBL}(b) we may apply Theorem \ref{TRescaledBBLWithWeights}, 
obtaining
\begin{eqnarray*}
	\mathcal L^{k+1}(Z_s(A,B)) &\geq & M^\frac{1}{k+1}_s \left(\left(\tilde{\tau}_{1-s}^{k, \alpha}(A,B)\right)^{k+1}\mathcal L^{k+1}(A),
	\left(\tilde{\tau}_{s}^{k, \alpha}(A,B)\right)^{k+1} \mathcal L^{k+1}(B)\right) \\
	&=& \left(\tau_{1-s}^{k,\alpha}(A,B)
	\mathcal L^{k+1}(A)^\frac{1}{k+1}+
	\tau_{s}^{k, \alpha}(A,B)
	\mathcal L^{k+1}(B)^\frac{1}{k+1}\right)^{k+1},
\end{eqnarray*}
which concludes the proof. {In a similar way, properties from (ii) and (iii) from Theorem \ref{TBrunn-Minkowski-2} can be obtained by Corollaries \ref{CLighterWeightedBBL} and \ref{CRescaledBBLWithoutWeights}, respectively.}
\end{remark}

\vspace{0.5cm} \noindent {\footnotesize{\sc  Mathematisches Institute,
		Universit\"at Bern,
		Sidlerstrasse 5,
		3012 Bern, Switzerland.}\\
	Email: \textsf{zoltan.balogh@math.unibe.ch}\\

	\noindent {\footnotesize {\sc Department of Economics, Babe\c s-Bolyai University, Str. Teodor Mihali 58-60, 400591
			Cluj-Napoca, Romania \& Institute of Applied Mathematics, \'Obuda University,
			B\'ecsi \'ut 96, 1034 Budapest, Hungary.}\\ Email:
		{\textsf{alex.kristaly@econ.ubbcluj.ro}}\\
		
		\noindent {\footnotesize{\sc   Mathematisches Institute,
				Universit\"at Bern,
				Sidlerstrasse 5,
				3012 Bern, Switzerland.\\
			} Email: {\textsf{kinga.sipos@math.unibe.ch}\\

\end{document}